\documentclass{amsart}
\usepackage{amssymb,amsmath, amsthm,latexsym}
\newcommand{\cal}[1]{\mathcal{#1}}
\theoremstyle{plain}

\parskip=\bigskipamount

\let\egthree=\phi
\let\phi=\varphi
\let\varphi=\egthree

\begin{document}
\title{Bounded cohomology and isometry groups of hyperbolic spaces}
\author{Ursula Hamenst\"adt}
\thanks
{{\it e-mail address:} ursula@math.uni-bonn.de\\
Partially supported by Sonderforschungsbereich 611}
\date{September 27, 2006}

\begin{abstract}
Let $X$ be an arbitrary hyperbolic geodesic metric space
and let $\Gamma$ be a countable
subgroup of the isometry group ${\rm Iso}(X)$ of $X$.
We show that if $\Gamma$ is non-elementary
and weakly
acylindrical (this is a weak
properness condition) then the second
bounded cohomology groups $H_b^2(\Gamma,\mathbb{R})$,
$H_b^2(\Gamma,\ell^p(\Gamma))$ $(1< p <\infty)$
are infinite dimensional.
Our result holds for example for
any subgroup of the mapping class group of
a non-exceptional surface of finite
type not containing
a normal subgroup which
virtually splits as a direct product.
\end{abstract}

\maketitle

\section{Introduction}

A \emph{Banach module} for a countable group
$\Gamma$ is a Banach space $E$ together with a
homomorphism of $\Gamma$ into the group of linear
isometries of $E$. For every such Banach module $E$
for $\Gamma$ and every $i\geq 1$, the
group $\Gamma$ naturally acts on the vector space
$L^\infty(\Gamma^i,E)$ of bounded functions $\Gamma^i\to E$.
If we denote by
$L^\infty(\Gamma^i,E)^\Gamma\subset L^\infty(\Gamma^i,E)$
the linear subspace of all $\Gamma$-invariant such functions, then the
\emph{second bounded
cohomology group} $H_b^2(\Gamma,E)$ of $\Gamma$
with coefficients $E$
is defined as the second cohomology group of the complex
\begin{equation} 0\to
L^\infty(\Gamma,E)^\Gamma \xrightarrow{d} L^\infty
(\Gamma^2,E)^\Gamma\xrightarrow{d} \dots
\end{equation} with the usual
homogeneous coboundary operator $d$ (see \cite{M}).
There is a natural
homomorphism of $H_b^2(\Gamma,E)$ into the ordinary
second cohomology group $H^2(\Gamma,E)$ of $\Gamma$
with coefficients $E$ which in general
is neither injective nor surjective.

In this paper we are only interested in the case that
$E=\mathbb{R}$ with the trivial $\Gamma$-action
or that $E=\ell^p(\Gamma)$
for some $p\in (1,\infty)$ with the
natural $\Gamma$-action by right translation which
assigns to a $p$-summable function $f$ and an element
$g\in \Gamma$ the function $gf:h\to f(hg)$.

Since every homomorphism $\rho$ of $\Gamma$ into a countable group
$G$ induces a homomorphism $\rho^*:H_b^2(G,\mathbb{R})\to
H_b^2(\Gamma,\mathbb{R})$, second bounded cohomology
with real coefficients can be used
to find obstructions to the existence of interesting homomorphisms
$\Gamma\to G$. The underlying idea is to find conditions
on $G$ and $\rho$ which ensure that the image
of the map $\rho^*$ is ``large'' (e.g. infinite dimensional)
and conclude that this imposes restrictions
on the group $\Gamma$.

Two countable groups
$\Gamma,G$ are called \emph{measure equivalent}
\cite{G93} if $\Gamma,G$ admit
commuting measure preserving actions
on a standard infinite measure Borel
space with finite measure fundamental domains.
Measure equivalence defines
an equivalence relation for countable groups \cite{Fu99a}.
Monod and Shalom
\cite{MS05} showed that for
countable groups, vanishing of the
second bounded cohomology groups with coefficients
in the regular representation is preserved under
measure equivalence. Thus
second bounded cohomology with coefficients in the regular
representation provides an obstruction
to the existence of a measure equivalence between
two given countable groups.

For the application of these ideas
it is necessary to obtain information on
these bounded cohomology groups. The first and easiest result
in this direction is due to B.~Johnson (see \cite{I} for
a discussion and references) who showed that
the bounded cohomology of
amenable groups with real coefficients is trivial.
Later Brooks \cite{B81} found a combinatorial
method for the construction of non-trivial real
second bounded cohomology
classes and used it to show that the second bounded
cohomology group of a finitely generated free
group is infinite dimensional.

Fujiwara \cite{F} investigated the second real
bounded cohomology group of a group of isometries of
a \emph{hyperbolic} geodesic metric space. Such
a space $X$ admits
a geometric
boundary $\partial X$. Each isometry of $X$ acts as a homeomorphism
on $\partial X$. The
\emph{limit set} of a group $\Gamma$ of isometries of $X$ is the
closed $\Gamma$-invariant subset of $\partial X$
of all accumulation points
of a fixed $\Gamma$-orbit in $X$. The group $\Gamma$
is called \emph{non-elementary} if its
limit set contains at least 3 points; then
the limit set of $\Gamma$ is in fact uncountable.
Using a refinement and an
extension of Brooks' method, Fujiwara showed that for
a countable non-elementary group $\Gamma$ of isometries of $X$
acting properly discontinuously on $X$ in a metric sense,
the kernel of the
map $H_b^2(\Gamma,\mathbb{R})\to H^2(\Gamma,\mathbb{R})$ is
infinite dimensional \cite{F}.
Bestvina and Fujiwara extended this
result further to countable subgroups of ${\rm Iso}(X)$
whose actions on $X$ satisfy some weaker properness assumption
\cite{BF}.
Their result is for example valid for
non-elementary subgroups of the
\emph{mapping class group} of an oriented surface $S$ of finite
type and negative Euler characteristic, i.e. for subgroups of the
group of isotopy classes of orientation preserving diffeomorphisms
of $S$ which are not virtually abelian. As a consequence, the
second bounded cohomology group of every non-elementary subgroup
of such a mapping class group is infinite dimensional.

On the other hand, by a result of Burger and Monod
\cite{BM99,BM02},
for every irreducible lattice $\Gamma$ in a connected
semi-simple Lie group with finite center, no compact factors
and of rank at least 2
the kernel of the natural map $H_b^2(\Gamma,\mathbb{R})\to
H^2(\Gamma,\mathbb{R})$ vanishes.
Together with the results of Fujiwara \cite{F} and
Bestvina and Fujiwara \cite{BF}
it follows easily that
the image of every
homomorphism of $\Gamma$ into a finitely generated
word hyperbolic group
or into the mapping class group of an oriented
surface of finite type
and negative Euler characteristic is finite \cite{BM02,BF}.
The latter result was earlier derived with different
methods by
Farb and Masur \cite{FM} building on the work of
Kaimanovich and Masur \cite{KM}.

The goal of this note is to present a new method for
constructing nontrivial second
bounded cohomology classes for a countable group $\Gamma$
from dynamical properties of suitable
actions of $\Gamma$.
We use it to give a common proof of
extensions of the above mentioned results of
Fujiwara \cite{F} and of Bestvina and Fujiwara
\cite{BF} which among other things
answers a question raised by Monod and Shalom
\cite{MS,MS05}.

For the formulation of these results, call a
countable group $\Gamma$ of isometries of
a (not necessarily proper)
hyperbolic geodesic metric
space $X$ \emph{weakly acylindrical} if for every point $x_0\in X$
and every $m>0$ there are numbers $R(x_0,m)>0,c(x_0,m)>0$ with the
following property. If $x,y\in X$ are such that a geodesic
$\gamma$ connecting $x$ to $y$ meets the $m$-neighborhood of $x_0$
and if $d(x,y)\geq R(x_0,m)$ then there are at most $c(x_0,m)$
elements $g\in \Gamma$ such that $d(x,gx)\leq m$ and $d(y,gy)\leq
m$ (compare with the definition of an
acylindrical isometry group in \cite{B03}).
We show in Section 4 (see \cite{F,BF,MMS} for
earlier results).

\bigskip

{\bf Theorem A:} {\it Let $\Gamma$ be a
non-elementary weakly acylindrical countable group
of isometries of an
arbitrary hyperbolic geodesic metric space. Then the kernels of the
maps $H_b^2(\Gamma,\mathbb{R})\to H^2(\Gamma,\mathbb{R})$ and
$H_b^2(\Gamma,\ell^p(\Gamma))\to H^2(\Gamma,\ell^p(\Gamma))$ $
(1< p<\infty)$ are infinite dimensional.}

\bigskip

As an easy corollary of Theorem A and a result of
Bowditch \cite{B03}
we obtain an extension of the
result of Bestvina and Fujiwara \cite{BF}. For its formulation,
we say that a group $\Gamma$ \emph{virtually splits} as a direct
product if $\Gamma$ has a finite index subgroup $\Gamma^\prime$
which splits as a direct product of two
infinite groups. We show.

\bigskip
{\bf Corollary B:} {\it Let $\Gamma$ be a subgroup of the mapping
class group of an oriented surface of finite type and negative
Euler characteristic. If $\Gamma$ is not virtually abelian
then the kernel of the map
$H_b^2(\Gamma,\mathbb{R})\to H^2(\Gamma,\mathbb{R})$
is infinite dimensional. If moreover $\Gamma$ does not
contain a normal subgroup which 
virtually split as a direct product then the kernel of
each of the maps
$H_b^2(\Gamma,\ell^p(\Gamma))\to H^2(\Gamma,\ell^p(\Gamma))$
$(1< p
<\infty)$ is infinite dimensional.}

\bigskip

The following corollary is
an immediate consequence of Corollary B and the
work of Burger-Monod and Monod-Shalom \cite{BM99,BM02,MS05}.
For its formulation, call a lattice $\Lambda$ in a product
$G=G_1\times G_2$ of two locally compact $\sigma$-compact
and non-compact topological groups \emph{irreducible}
if the projection of $\Lambda$ to each of the factors is dense.

\bigskip

{\bf Corollary C:} {\it Let $\Gamma$ be a
subgroup of the mapping class group of
an oriented surface of finite type and negative
Euler characteristic. Assume that
$\Gamma$ does not contain a normal subgroup which
virtually splits as a direct product.
Then $\Gamma$ is not measure equivalent to an
irreducible lattice in a product of two
locally compact $\sigma$-compact non-compact topological
groups.}

\bigskip

For lattices in semisimple Lie
groups of higher rank, Corollary C follows from   
\cite{FM} and the beautiful work of
Furman \cite{Fu99a}.
The earlier result of
Zimmer \cite{Z91} suffices to deduce Corollary C
for the full mapping class group which admits
a linear representation with infinite image.
Recently, Kida \cite{K06}
derived a much stronger rigidity result. Namely,
he showed that for every countable group $\Lambda$ which
is measure equivalent to the mapping class
group ${\cal M}$ 
of a non-exceptional oriented surface of finite type,
there is a homomorphism $\Lambda\to {\cal M}$ with
finite kernel and finite index image.

The organization of the paper is as follows.
In Section 2, we introduce
our method for the construction of
second bounded cohomology classes in the concrete example of
the fundamental group $\Gamma$ of a \emph{convex cocompact}
complete Riemannian manifold $M$ of bounded negative
curvature. Such a manifold $M$ contains a
compact convex subset ${\cal C}(M)$, the so-called \emph{convex
core}, as a strong deformation retract. The group
$\Gamma$ is a word hyperbolic, and the convex core
${\cal C}(M)$ of $M$ is a
$K(\Gamma,1)$-space. Therefore,
if $\Gamma$ is \emph{non-elementary},
i.e. if $\Gamma$ is not abelian,
then the dimension of
the cohomology group $H^2(\Gamma,\mathbb{R})$ is finite, and
by Fujiwara's result \cite{F}, the
group $H_b^2(\Gamma,\mathbb{R})$ is infinite dimensional.

Inspired by
a result of Barge and Ghys \cite{BG}, we
relate the second bounded cohomology group $H_b^2(\Gamma,\mathbb{R})$
to the \emph{geodesic flow}
$\Phi^t$ of $M$ which
acts on the unit tangent bundle $T^1M$ of $M$. Since $M$
is convex cocompact, $T^1M$
admits a compact $\Phi^t$-invariant hyperbolic
set $W$ which is the closure of
the union of all closed orbits of $\Phi^t$.
The restriction of $\Phi^t$ to $W$ is topologically transitive.

A \emph{cocycle} for the action of $\Phi^t$ on $W$ is a continuous
function $c:W\times \mathbb{R}\to \mathbb{R}$ such that
$c(v,s+t)=c(v,t)+c(\Phi^tv,s)$ for all $v\in W$ and all $s,t\in
\mathbb{R}$. Two cocycles $c,d$ are \emph{cohomologous} if there
is a continuous function $\psi:W\to \mathbb{R}$ such that
$\psi(\Phi^tv)+c(v,t)= d(v,t)+\psi(v)$. The collection of all
cocycles which are cohomologous to a given cocycle $c$ is the
\emph{cohomology class} of $c$. The \emph{flip} ${\cal F}:v\to -v$
acts on $W$ and on the space of cocycles for the geodesic flow
preserving cohomology classes. The cohomology class of a cocycle
$c$ is called \emph{flip anti-invariant} if ${\cal F}(c)$ is
cohomologous to $-c$.
We denote by ${\cal D\cal C}(M)$ the vector
space of all flip anti-invariant cohomology classes of
cocycles for the geodesic flow on $W$ which are H\"older
continuous, i.e. such that for a fixed number $t>0$ the
function $v\to c(v,t)$ is H\"older continuous.

Every smooth
closed 1-form on $M$ defines via integration along
orbit segments of $\Phi^t$
a H\"older continuous
cocycle for $\Phi^t$ which is anti-invariant under
the flip. Two cocycles defined by closed one-forms in this
way are cohomologous if and only if the one-forms
define the same de Rham cohomology class on $M$. Thus
$H^1(\Gamma,\mathbb{R})=H^1(M,\mathbb{R})$ is naturally a subspace
of ${\cal D\cal C}(M)$. We show in Section 2.

\bigskip

{\bf Theorem D:} {\it Let $\Gamma$ be the fundamental group of a
convex cocompact manifold $M$ of bounded negative curvature. Then
the quotient space ${\cal D\cal C}(M)/H^1(\Gamma,\mathbb{R})$
naturally embeds into ${\rm ker}(H_b^2(\Gamma,\mathbb{R})\to
H^2(\Gamma,\mathbb{R}))$.}

\bigskip

Section 3 contains the main technical result of this paper.
Embarking from the concrete construction in Section 2, we present
an abstract dynamical criterion for infinite dimensional second
bounded cohomology  for a countable group $\Gamma$ acting as a
group of homeomorphisms on a metric space of bounded diameter. The
coefficients of these cohomology groups can be either $\mathbb{R}$
or $\ell^p(\Gamma)$ for some $p\in (1,\infty)$. Theorem 4.4 of
Section 4 shows that our criterion can be applied to countable
groups which admit a non-elementary weakly acylindrical isometric
action on a hyperbolic geodesic metric space; this then yields
Theorem A. Section 5 contains
the proof of Corollary B and Corollary C as well as
a short discussion of some applications which are due
to Monod and Shalom.

\section{Dynamical cocycles and bounded cohomology}

In this section we consider an $n$-dimensional
convex cocompact Riemannian
manifold $M$ of bounded negative curvature.
Then $M=\tilde M/\Gamma$ where
$\tilde M$ is a simply connected complete Riemannian manifold
of bounded negative curvature and
$\Gamma$ is a group
of isometries acting properly discontinuously and freely
on $\tilde M$. The manifold $\tilde M$ admits a natural
compactification by adding the \emph{geometric boundary}
$\partial \tilde M$ which is a topological sphere of
dimension $n-1$. Every isometry of $\tilde M$ acts on
$\partial \tilde M$ as a homeomorphism.
The \emph{limit set} $\Lambda$ of $\Gamma$ is
the set of accumulation points in $\partial
\tilde M$ of a fixed
orbit $\Gamma x$ $(x\in \tilde M)$
of the action of $\Gamma$ on $\tilde M$.
We always assume that the group $\Gamma$ is non-elementary,
i.e. that its limit set contains as least 3 points.
Then $\Lambda$ is
the smallest nontrivial closed subset of $\partial\tilde M$ which is
invariant under the action of $\Gamma$.

The \emph{geodesic flow} $\Phi^t$ acts on the
\emph{unit tangent
bundle} $T^1\tilde M$ of $\tilde M$ and on the
unit tangent bundle $T^1M$ of $M$.
Let $\tilde L\subset T^1\tilde M$ be the set of all
unit tangents of all geodesics with both end-points in
$\Lambda$. Then $\tilde L$ is invariant under the
action of $\Phi^t$ and the action of $\Gamma$.
The quotient $L=\tilde L/\Gamma$ is just the
\emph{non-wandering set} for the action of $\Phi^t$
on $T^1 M$.
Since $M$ is convex cocompact by
assumption, $L$ is a compact
hyperbolic set for the geodesic flow $\Phi^t$ on
$T^1M$. The sets $\tilde L$ and $L$ are moreover
invariant under the \emph{flip} ${\cal F}:v\to -v$ which maps a
unit tangent vector to its negative.
The Riemannian metric on $M$ induces a complete Riemannian metric
and hence a complete distance function $d$ on $T^1M$.

A continuous real-valued \emph{cocycle}
for the action of $\Phi^t$ on $L$ is
a continuous function $c:L\times \mathbb{R}\to \mathbb{R}$ with
the property that $c(v,t+s)=c(v,t)+c(\Phi^t v,s)$ for all $v\in L$,
all $s,t\in \mathbb{R}$. Every continuous function
$f:L\to \mathbb{R}$ defines such a
cocycle $c_f$ by
$c_f(v,t)=\int_0^t f(\Phi^s v)ds$. Two
cocycles $b,c$ are called \emph{cohomologous} if there is a
continuous function $\psi:L\to \mathbb{R}$ such that $\psi(\Phi^t
v)+c(v,t)-\psi(v)=b(v,t)$.
If $b, c$ are
cocycles which are H\"older
continuous with respect to
the distance $d$ on $L$, i.e. if for fixed $t>0$ the maps
$v\to b(v,t),v\to c(v,t)$ are H\"older continuous,
then by Livshicz' theorem, $b,c$ are
cohomologous if and only if for every periodic point $v$ of the
geodesic flow with period $\tau >0$ we have $b(v,\tau)=c(v,\tau)$
\cite{HK}. Every H\"older
continuous cocycle is cohomologous to the cocycle of
a H\"older continuous function
$f$ (see e.g. \cite{H99}), and
two H\"older functions $f,g$ on $L$ are
\emph{cohomologous}, i.e. their cocycles $c_f,c_g$
are cohomologous, if and
only if we have $\int_{\gamma^\prime} f=\int_{\gamma^\prime}g$ for
every closed geodesic $\gamma$ on $M$ (where $\gamma^\prime$ is the
unit tangent field of $\gamma$).

The flip ${\cal F}$ acts on the space of cocycles preserving
cohomology classes. We call the cohomology class of a cocycle $c$
\emph{anti-invariant} under the flip ${\cal F}$ if ${\cal F}(c)$
is cohomologous to $-c$. If the cohomology class of a H\"older
continuous cocycle $c$ is anti-invariant under the flip then there
is a H\"older continuous function $f$ which is anti-invariant under
the flip, i.e. which satisfies
$f(v)=-f(-v)$ for all $v\in L$, such
that the cocycle $c_f$ defined by $f$ is cohomologous to $c$
(compare \cite{H97}). Denote by ${\cal A}$ the vector space of all
H\"older continuous functions $f$ on $L$ which are anti-invariant
under the flip ${\cal F}$.

Since $L$ is a compact invariant hyperbolic
topologically transitive
set for the geodesic flow on
$T^1M$, for every H\"older continuous function $f$ on $L$ and
every number $\epsilon_0 >0$ which is smaller than half of the
injectivity radius of $M$ there is a number $k>0$ only depending
on the H\"older norm of $f$ with the following property.
Let $v,w\in L$ and let $T>0$ be
such that
$d(\Phi^tv,\Phi^t w)\leq \epsilon_0$ for all $t\in
[0,T]$; then
$\vert \int_0^Tf(\Phi^t
v)dt-\int_0^Tf(\Phi^t w)dt\vert \leq k$.

A \emph{quasi-morphism} for $\Gamma$ is a function
$\phi:\Gamma\to \mathbb{R}$
such that
\begin{equation}
\Vert \phi\Vert_0=
\sup_{g,h\in \Gamma}\vert \phi(g)+\phi(h)-\phi(gh)\vert<\infty.
\end{equation}
The set ${\cal Q}$ of
all quasi-morphisms for $\Gamma$ naturally has the structure of a
vector space. The function $\Vert \,\Vert_0:{\cal Q}\to
[0,\infty)$ which associates to a quasi-morphism $\phi$
its \emph{defect} $\Vert \phi\Vert_0$ is a pseudo-norm
which vanishes precisely on the subspace of \emph{morphisms}.

\bigskip

{\bf Lemma 2.1:} {\it There is a linear map $\Psi:{\cal A}\to
{\cal Q}$. For every $f\in {\cal A}$, the defect
$\Vert \Psi(f)\Vert_0$ of $\Psi(f)$ is bounded from
above by a constant only depending on the curvature
bounds of $M$ and 
the H\"older norm of $f$.}

{\it Proof:} Let $f\in {\cal A}$ and extend $f$ to a
locally H\"older continuous
flip anti-invariant function $F$ on $T^1M$.
Such an extension can be
constructed as follows.
Choose a probability measure
$\mu$ on $L$ for which there are constants
$0<a<b$ such that the $\mu$-mass of a ball
$B(v,r)$ of radius $r<1$ about
a point $v\in L$ is contained in $[r^b,r^a]$; for example, the
unique measure of maximal entropy for the geodesic flow on $L$
has this property.
We view $\mu$ as a probability measure on $T^1M$ which is
supported in $L$. Let $\tau:[0,\infty)\to
[0,1]$ be a smooth function which satisfies
$\tau(t)=1$ for $t$ close to $0$ and $\tau[1,\infty)=0$.
Via multiplying the restriction of
$\mu$ to the metric ball $B(w,r)$ $(w\in T^1M)$ by
the function $z\to \tau(d(z,w)/r)$
we may assume that the measures
$\mu\vert B(w,r)$ depend continuously on $w\in T^1M,r>0$ in
the weak$^*$-topology.

For $w\in T^1M$ let $\delta(w)\geq 0$ be the distance between $w$
and $L$. For $w\in T^1M -L$ define
\begin{equation}
f_0(w)=\int_{B(w,2\delta(w))\cap L}f d\mu/\mu(B(w,2\delta(w)))
\end{equation}
and let $f_0(w)=f(w)$ for $w\in L$.
By assumption on
the measures $\mu\vert B(w,r)$ and since
$f$ is H\"older continuous, the
function $f_0$ is locally H\"older continuous
and its restriction to $L$
coincides with $f$. Thus we obtain a locally H\"older continuous flip
anti-invariant extension $F$ of $f$ to $T^1M$ by assigning to
$w\in T^1M-L$ the value $F(w)=\frac{1}{2}(f_0(w)-f_0(-w))$.
For every compact subset $K$ of $T^1M$ the H\"older norm
of the restriction of $F$ to $K$ only depends on $K$ and
on the H\"older norm of $f$.
If
$F,H$ are the extensions of $f,h$ constructed in
this way and if $a,b\in \mathbb{R}$ then $aF+bH$ is the extension
of $af+bh$.

Let again $\Lambda$ be the limit set of $\Gamma$.
The closure
${\rm Conv}(\Lambda)\subset \tilde M$ of the
convex hull of $\Lambda$ in $\tilde M$
is invariant under
the action of $\Gamma$. The \emph{convex core}
${\cal C}(M)={\rm Conv}(\Lambda)/\Gamma$
of $M$ is compact.
Let $\tilde F$ be the lift of $F$ to $T^1\tilde M$ and choose a point
$p\in {\rm Conv}(\Lambda)$.
For an element $g\in \Gamma$ define $\Psi(f)(g)$ to be
the integral of $\tilde F$ over the tangent of the oriented
geodesic joining $p$ to $g(p)$. We claim that $\Psi(f)$ is a
quasi-morphism for $\Gamma$.

For the proof of this claim,
recall that the curvature of $\tilde M$ is pinched
between two negative constants and therefore
by comparison,
for every $\epsilon_0>0$ there is
a number $k=k(\epsilon_0)>0$ only depending on $\epsilon_0$ and an
upper curvature bound for $\tilde M$ with the following property.
Let $T$ be a geodesic triangle in $\tilde M$ with
vertices $A_1,A_2,A_3$ and denote by $a_i$ the side of $T$
connecting $A_{i-1}$ to $A_{i+1}$. Let
$q_i\in a_i$ be the nearest point
projection of the vertex $A_{i}$ to the side $a_i$ and let
$\gamma_{i,+},\gamma_{i,-}$ be the
oriented geodesic
arc parametrized by arc length which
connects $q_i=\gamma_{i,+}(0)$ to $A_{i+1}=\gamma_{i,+}(\tau_{i,+})$,
$q_i=\gamma_{i,-}(0)$ to $A_{i-1}=\gamma_{i,-}(\tau_{i,-})$
(here indices are taken modulo 3).
Then $t_i=\tau_{i+1,+}-\tau_{i,-}\in [-k,k]$ and moreover
for every $t\in [k,\tau_{i,-}]$ the distance between
$\gamma_{i,-}^\prime(t)\in T^1\tilde M$  and
$\gamma_{i+1,+}^\prime(t+t_i)\in T^1\tilde M$
is at most $\epsilon_0$.

Now by assumption, the function $F$ is anti-invariant under the
flip and locally H\"older continuous.
Therefore our above discussion implies that
the integral of $\tilde
F$ over the unit tangent field of a closed
curve in ${\rm Conv}(\Lambda)$
which consists of three geodesic arcs forming a geodesic
triangle is bounded from above in
absolute value by a universal constant times the H\"older norm
of the restriction of $F$ to the compact subset
$T^1M\vert {\cal C}(M)$ of $T^1M$ of all unit
vectors with foot point in
the convex core ${\cal C}(M)={\rm Conv}(\Lambda)/\Gamma$.
On the other hand, by invariance of $\tilde F$ under the
action of $\Gamma$ and by anti-invariance of
$\tilde F$ under the flip, for $g,h\in \Gamma$ the
quantity $\vert \Psi(f)(g)+\Psi(f)(h)-\Psi(f)(gh)\vert$ is just
the absolute value of the integral of $\tilde F$ over
the unit tangent field of the oriented
geodesic triangle
in ${\rm Conv}(\Lambda)\subset \tilde M$
with vertices $p,g(p),g(h (p))$.
Thus $\Psi(f)$ is indeed a quasi-morphism and the assignment $f\to
\Psi(f)$ defines a linear map $\Psi:{\cal A}\to {\cal Q}$. 
Moreover, the defect $\Vert \phi\Vert_0$ of $\phi$ is bounded
from above by a constant only depending on the 
curvature bounds of $M$ and the H\"older norm 
of $f$. This
shows the lemma. \qed

\bigskip

Two quasi-morphisms $\phi,\psi$ for $\Gamma$ are called
\emph{equivalent} if $\phi-\psi$ is a bounded function. This is
clearly an equivalence relation.
If $\phi_1$ is equivalent to $\phi_2$ and $\psi_1$ is equivalent
to $\psi_2$ then for all $a,b\in \mathbb{R}$ the quasi-morphism
$a\phi_1+b\psi_1$ is equivalent to
$a\phi_2+b\psi_2$ and hence
the set ${\cal Q}B$ of
equivalence classes of quasi-morphisms of $\Gamma$
has a natural structure of a vector space. It
contains as a
subspace the vector space $H^1(\Gamma,\mathbb{R})$ of all equivalence
classes of \emph{morphisms} of $\Gamma$.
There is an exact sequence
\begin{equation}
0\to H^1(\Gamma,\mathbb{R})\to {\cal Q}B\to
H_b^2(\Gamma,\mathbb{R}) \to H^2(\Gamma,\mathbb{R})
\end{equation}
and therefore the quotient space $\tilde {\cal Q}={\cal Q}B
/H^1(\Gamma,\mathbb{R})$ can naturally be identified with the
kernel of the map $H_b^2(\Gamma,\mathbb{R})\to
H^2(\Gamma,\mathbb{R})$ (see \cite{M}).
In particular, an equivalence class of
quasi-morphisms can be viewed as a cohomology class of
$\Gamma$-invariant bounded cocycles $\phi\in
L^\infty(\Gamma^3,\mathbb{R})^\Gamma$. In this interpretation,
the cocycle $\phi$ determined by the quasi-morphism
$\psi$ associates to a triple $(g,h,u)\in \Gamma^3$
the value $\phi(g,h,u)=\psi(g^{-1}h)+\psi(h^{-1}u)-\psi(g^{-1}u)$.

For $f\in {\cal A}$, the definition of the quasi-morphism
$\Psi(f)$ in the proof of Lemma 2.1
depends on the choice of an extension of $f$ to a
locally H\"older continuous
flip anti-invariant function
on $T^1M$ and also on the choice of a
basepoint $p\in {\rm Conv}(\Lambda)$. The next lemma shows that
the cohomology class of $\Psi(f)$ is independent of these choices.

\bigskip

{\bf Lemma 2.2:} {\it The cohomology class of the quasi-morphism
$\Psi(f)$ does not depend on the choice of the basepoint $p$ nor
on the extension of $f$ to a locally H\"older continuous
flip anti-invariant function on $T^1M$.}

{\it Proof:} Let $f\in {\cal A}$ and let $F$ be a
locally H\"older
continuous flip anti-invariant
extension of $f$ to $T^1M$. Denote by $\Psi(f)$ the
quasi-morphism constructed in the proof of Lemma 2.1 using the
extension $F$ of $f$ and the basepoint $p\in {\rm Conv}(\Lambda)$.
We first show that a different choice $q\in {\rm Conv}(\Lambda)$
of a basepoint gives rise to a quasi-morphism which is equivalent
to $\Psi(f)$ in our above sense. For this we follow
\cite{BG}. Let
$\tilde F$ be the lift of $F$ to $T^1\tilde M$. For $g\in \Gamma$
define $\rho(g)$ to be the integral of the function
$\tilde F$ over the unit tangent field of the oriented geodesic
quadrangle with vertices
$q,g(q),g(p),p$. As in the proof of Lemma 2.1 we conclude that the
function $\rho:\Gamma\to \mathbb{R}$ is uniformly bounded. Since
by invariance of $\tilde F$ under the action of $\Gamma$
the integral of $\tilde F$ over the oriented geodesic arc connecting
$g(q)$ to $g(p)$ is independent of $g\in
\Gamma$ and, in particular, it coincides with the negative of the integral
of $\tilde F$ over the oriented geodesic arc
connecting $p$ to $q$, the
quasi-morphism defined by $F$ and the basepoint $q$
just equals $\Psi(f)+\rho$. Thus changing
the basepoint does not change the equivalence class of our
quasi-morphism $\Psi(f)$.

Now we may also replace the point $q\in {\rm Conv}(\Lambda)$ by a
point $\xi\in \Lambda$. Namely, for $g\in \Gamma$ define $\rho(g)$
to be the oriented integral of the function $\tilde F$ over the
tangent lines of the ideal geodesic quadrangle with vertices
$\xi,g(\xi),g(p),p$. As before, this function is uniformly
bounded. By invariance of $\tilde F$
under the action of $\Gamma$ and the fact
that $\tilde F$ is anti-invariant under the flip
we obtain that the 2-cocycle
for $\Gamma$ defined as above by the quasi-morphism $\Psi(f)+\rho$ is
just the cocycle $\eta\in
L^\infty(\Gamma^3,\mathbb{R})^\Gamma$ which assigns to a triple
$(g,h,u)\in \Gamma^3$ the integral of $\tilde F$ over the
unit tangents of the (possibly degenerate) oriented
ideal triangle with vertices $g(\xi), h(\xi),u(\xi)$.
Since these unit tangents are contained in the lift
$\tilde L$ of the non-wandering set
$L$ for the geodesic flow on $T^1M$, the
cocycle $\eta$ only depends on $f$ but not on an extension of $f$
to $T^1M$. Thus the cohomology class defined by $\Psi(f)$ is
independent of the extension as well. \qed

\bigskip

In the sequel we denote for $f\in {\cal A}$
by $\Theta(f)\in \tilde{\cal Q}\sim{\rm
ker}(H_b^2(\Gamma,\mathbb{R})\to H^2(\Gamma,\mathbb{R}))$ the
cohomology class of the quasi-morphism $\Psi(f)$. By Lemma 2.2,
this class only depends on $f$.
Moreover, the assignment $\Theta:{\cal A}\to \tilde Q$ is clearly
linear. We next investigate the kernel of the map $\Theta$.

Since $\Gamma$ is convex cocompact by assumption,
there is a natural correspondence between
oriented closed geodesics on $M$
and conjugacy classes in $\Gamma$. For every homomorphism
$\rho:\Gamma\to \mathbb{R}$ and every $g\in \Gamma$,
the value $\rho(g)$ of $\rho$ on $g$
only depends on the conjugacy class of $g$. Therefore such
a homomorphism defines a function on the set of closed geodesics
on $M$; we denote the value of $\rho$ on such a closed geodesic
$\gamma$ by $\rho(\gamma)$. We have.

\bigskip

{\bf Lemma 2.3:} {\it $\Theta(f)=0$ if and only if there is a
morphism $\rho:\Gamma\to \mathbb{R}$ such that
$\int_{\gamma^\prime} f =\rho(\gamma)$ for every closed geodesic
$\gamma$ on $M$.}

{\it Proof:} Let $f\in {\cal A}$ and assume that there is a
morphism $\rho:\Gamma\to \mathbb{R}$ such that
$\int_{\gamma^\prime} f=\rho(\gamma)$ for every closed geodesic
$\gamma$ on $M$. This morphism defines a class in
$H^1(M,\mathbb{R})$ and therefore by the de Rham theorem,
there is a smooth closed 1-form $\beta$ on $M$ which
defines $\rho$ via integration.
Let $\tilde \beta$ be the pull-back of $\beta$ to a closed
$1$-form on $\tilde M$. Then $\tilde \beta$ is exact and hence
the integral of $\tilde \beta$
over every piecewise smooth closed curve
in $\tilde M$ vanishes.

By Livshicz' theorem \cite{HK} and the choice of $\beta$,
there is a H\"older continuous
flip anti-invariant function $\psi:L\to \mathbb{R}$ such that
$\int_0^Tf(\Phi^t v) dt= \psi(\Phi^T v)+\int_0^T\beta(\Phi^t
v)dt-\psi(v)$ for every $v\in L$ and all $T>0$.
As in the proof of Lemma 2.1 we extend $\psi$ to a
locally H\"older continuous 
flip anti-invariant function on all of $T^1M$
which we denote by the same symbol.
Let $\tilde \psi$ be the lift of $\psi$ to $T^1\tilde M$. Fix
a point $p\in {\rm Conv}(\Lambda)$ and for $g\in \Gamma$ let
$\gamma_g$ be the geodesic arc connecting $p=\gamma_p(0)$ to
$g(p)=\gamma_p(T)$. Define
$\alpha(g)=\tilde\psi(\gamma_g^\prime(T))+\int_0^T
\tilde\beta(\gamma_g^\prime(t))dt
-\tilde \psi(\gamma_g^\prime(0))$. By Lemma 2.1 and Lemma 2.2, $\alpha$
is a quasi-morphism for $\Gamma$ which defines the cohomology class
$\Theta(f)$. On the other hand, $\alpha$ differs from the
quasi-morphism defined by $\beta$ by a bounded function.
Since the integral of $\tilde\beta$ over every
piecewise smooth closed curve in
$\tilde M$ vanishes,
the cohomology class
$\Theta(f)$ of the quasi-morphism $\alpha$ vanishes.

On the other hand, let $f\in {\cal A}$ and assume that there is no
morphism $\rho:\Gamma\to \mathbb{R}$ such that
$\int_{\gamma^\prime}f=\rho(\gamma)$ for every closed geodesic
$\gamma$ on $M$. We have to show that the cohomology class
$\Theta(f)$
does not vanish. Using the exact sequence (4)
above, this is the case if and only if a quasi-morphism
$\Psi(f)$ representing $\Theta(f)$
is not equivalent to any morphism for $\Gamma$.

For this we argue as before. Namely,
let $\rho:\Gamma\to \mathbb{R}$ be any morphism
for $\Gamma$ and let $\beta$ be a smooth
closed 1-form on $M$ defining
$\rho$. By assumption, there is
a periodic point $v\in L$
of period $\tau>0$
for the geodesic flow
$\Phi^t$ such that $\int_0^\tau
(f-\beta)(\Phi^tv)dt=c>0$. Let $\tilde v$ be a lift of $v$ to
$\tilde L$ and let $p\in \tilde M$ be the foot-point of $\tilde
v$. Choose an extension of $f$ to a
locally H\"older continuous function
$F$ on $T^1M$ and let $\tilde F$ be the lift of $F$ to $T^1\tilde
M$. By definition, the quasi-morphism $\Psi(f)$ induced by $F$
and the choice of the basepoint $p$ assigns to $g\in \Gamma$
the integral $\int_0^T \tilde F(\gamma_g^\prime(s))ds$
where $\gamma_g:[0,T]\to
\tilde M$ is the oriented geodesic arc connecting $p$ to $g(p)$.
Moreover, this quasi-morphism represents the class $\Theta(f)$.
Now let $\eta$ be the geodesic in $\tilde M$ which is tangent to
$\tilde v$. By the choice of $\tilde v$ there is an element $h\in
\Gamma$ which preserves $\eta$ and whose restriction to $\eta$
is the translation $\eta(t)\to \eta(t+\tau)$
with translation length $\tau$. Hence we have
$\Psi(f)(h^m)=m\int_0^\tau f(\Phi^t v)dt$ and
$(\Psi(f)-\rho)(h^m)= mc$ for all $m\in \mathbb{Z}$.
In particular, the function $\Psi(f)-\rho$ is
unbounded and consequently $\Psi(f)$ is not equivalent to $\rho$. Since
$\rho$ was arbitrary this means that the projection of $\Psi(f)$
into $\tilde {\cal Q}$ does not vanish. \qed

\bigskip

Fix a number $\epsilon_0>0$ which is smaller than half of the
injectivity radius of $M$ and for $f\in {\cal A}$ define $\Vert
f\Vert_{\cal A}$ to be the infimum of the numbers $k>0$ with the property
that $\vert \int_0^T f(\Phi^t v)dt-\int_0^T f(\Phi^t w)dt \vert
\leq k$ whenever $v,w\in L$ and $T>0$ are such that
$d(\Phi^t v,\Phi^t w)\leq \epsilon_0$ for every
$t\in [0,T]$. We have.

\bigskip

{\bf Lemma 2.4:} {\it $\Vert \,\Vert_{\cal A}$ is a norm on ${\cal A}$.}

{\it Proof:} We observed above that $\Vert f\Vert_{\cal A}<\infty$ for
every H\"older continuous function $f\in {\cal A}$. Moreover, we
clearly have $\Vert a f\Vert_{\cal A}= \vert a\vert \Vert f\Vert_{\cal A}$ for
all $f\in {\cal A}$ and all $a\in \mathbb{R}$ and $\Vert
f+g\Vert_{\cal A}\leq \Vert f\Vert_{\cal A} +\Vert g\Vert_{\cal A} $ 
by a simple
application of the triangle inequality. Thus we are left with
showing that $\Vert f\Vert_{\cal A} =0$ only if $f\equiv 0$. For this
assume that $0\not\equiv f\in {\cal A}$. Since $f$ is
anti-invariant under the flip by assumption, $f$ is not the
constant function. Hence by continuity there are points $v,w\in W$
with $d(v,w)<\epsilon_0/2$ and numbers $\delta >0$, $T\in
(0,\epsilon_0)/2$ with $f(\Phi^tv)\geq f(\Phi^tw)+\delta$ for all
$t\in [0,T]$. Then $\Vert f\Vert_{\cal A} \geq \delta T$ by the
definition of $\Vert\,\Vert_{\cal A}$. \qed

\bigskip

Call two H\"older functions $f,g\in {\cal A}$ \emph{weakly
cohomologous} if $f-g$ is cohomologous to
a closed one-form
on $M$, viewed as a function on $T^1M$.
The class of $f$ under the thus defined equivalence
relation will be called the \emph{weak cohomology class} of $f$.
The set ${\cal H}$ of weak cohomology classes of H\"older
functions is a vector space.
For $\psi\in {\cal H}$ let
$\Vert \psi\Vert $ be the infimum of the norms $\Vert f\Vert
_{\cal A}$ where $f$ runs through all functions in ${\cal A}$ which
define the weak cohomology class $\psi$. Then $\Vert \,\Vert$ is a
pseudo-norm on ${\cal H}$.

The \emph{Gromov norm} $\Vert\,\Vert$ of an element
$\alpha\in H_b^2(\Gamma, \mathbb{R})$ is the infimum of the
supremums-norms over all bounded 2-cocycles for $\Gamma$
representing $\alpha$ \cite{G83} (here a bounded 2-cocycle
is a bounded $\Gamma$-invariant function on
$\Gamma^3$ contained in the kernel of the coboundary operator).
If $\phi:\Gamma\to \mathbb{R}$
is a quasi-morphism
then the Gromov norm of the cohomology
class defined by $\phi$ is the infimum of the defects
$\Vert\eta\Vert_0$ where $\eta$ runs through the collection of all
quasi-morphisms with the property that $\eta-\phi$ is
equivalent to a morphisms of $\Gamma$.

By Lemma 2.3, the map $\Theta$ factors through an injective linear
map ${\cal H} \to \tilde Q={\rm ker}(H_b^2(\Gamma,\mathbb{R})\to
H^2(\Gamma,\mathbb{R}))$ which we denote again by $\Theta$. The
following corollary summarizes our discussion and implies Theorem
D from the introduction.

\bigskip

{\bf Corollary 2.5:} {\it The map $\Theta:({\cal H},
\Vert\,\Vert)\to (\tilde{\cal Q}, \Vert\,\Vert)$ is a continuous
embedding.}

{\it Proof:} By Lemma 2.1, Lemma 2.3 and Lemma
2.4 we only have to show continuity of
$\Theta$.
For this choose a point $\xi\in \Lambda$.
Let $f\in {\cal A}$ and let $\tilde f$ be
the lift of $f$ to $T^1\tilde M$.
For $g,h,u\in \Gamma$
define $\alpha(g,h,u)$ to be the integral of
$\tilde f$ over the unit tangents of the (possibly
degenerate) oriented ideal
triangle with vertices $g\xi,h\xi,u\xi$.
By Lemma 2.3 and its proof,
$\alpha$ is a
cocycle which represents the class $\Theta(f)$.
The considerations
in the proof of Lemma 2.1 show that
$\vert\alpha(g,h,u)\vert\leq
c\Vert f \Vert_{\cal A}$ for a universal constant $c>0$,
in particular we have $\alpha\in L^\infty(\Gamma^3,\mathbb{R})^\Gamma$
and and the Gromov norm of the cohomology class defined
by $\alpha$ is not bigger than
$c\Vert f\Vert_{\cal A}$. From this
continuity of the map $\Theta$ follows. \qed

\section{A dynamical criterion for infinite dimensional
second bounded cohomology}

This section contains the main technical result of this note. We
consider an arbitrary countable group $\Gamma$ which acts by
homeomorphisms on a metric space $(X,d)$ of finite
diameter without isolated points.
Our goal
is to construct bounded cohomology classes for $\Gamma$ using
dynamical properties of the action of $\Gamma$ on $X$
as in Section 2. In the application we have in mind, the
space $X$ is the Gromov boundary of a hyperbolic geodesic
metric space and $\Gamma$ is a group of isometries 
acting on $X$ as a group of homeomorphisms.

We begin with
describing some properness condition for the
action of a countable group $\Gamma$ by homeomorphisms
of $(X,d)$. Namely,
the metric $d$ on our space $X$ induces a metric on the space
$X^3$ of triples of points in $X$ which we denote
again by $d$; this metric is given by
$d((x_1,x_2,x_3),(y_1,y_2,y_3))=\sum_{i=1}^3 d(x_i,y_i)$.
Let $\Delta\subset X^3$ be the closed subset
of all triples for which
at least two points in the triple coincide.
The diagonal
action of $\Gamma$ on $X^3$
preserves the open set $X^3-\Delta$ of triples of
pairwise distinct points in $X$.

{\bf Definition:} The action of $\Gamma$ on $X^3-\Delta$ is called
\emph{metrically proper} if for every $\nu\in (0,\frac{1}{2})$
there are constants $m(\nu)>0$, $R(\nu)>-\log \nu/4$
such that for any two open disjoint sets
$U,V\subset X$ of distance at least $\nu$ and of
diameter at most $e^{-R(\nu)}$ the following is satisfied.
\begin{enumerate}
\item Let $W\subset X$ be a set of diameter
at most $e^{-R(\nu)}$ whose distance to $U\cup V$ is
at least $\nu$. Write
$C=U\times V\times W\subset X^3-\Delta$; then
for all $k\in \mathbb{Z}$ and every fixed pair of points
$x_0\not=y_0\in X$ with $d(x_0,y_0)\geq \nu$
there are at most $m(\nu)$
elements $g\in \Gamma$ with
\begin{align}
g(C)\cap & \{(x,y,z)\in X^3-\Delta\mid x=x_0,y=y_0,
\notag \\
e^{-k} \leq & \min\{d(z,x_0),d(z,y_0)\}\leq
e^{-k+1} \}\not=\emptyset. \notag
\end{align}
\item Let $U^\prime,V^\prime\subset X$ be open disjoint
sets of distance at least $\nu$ and of diameter at most
$e^{-R(\nu)}$. Let
$Z\subset X$ (or $Z^\prime \subset X$) be the set of all points
whose distance to $U\cup V$ (or to $U^\prime\cup V^\prime$)
is bigger than $\nu$.
Then there are at most
$m(\nu)$ elements $g\in \Gamma$ with
\begin{equation}
g(U\times V\times Z)\cap U^\prime\times V^\prime\times
Z^\prime\not=\emptyset.
\notag
\end{equation}
\end{enumerate}

If the action of $\Gamma$ on $X^3-\Delta$ is
metrically proper, then every point in $X^3-\Delta$ has a
neighborhood $N$ in $X^3-\Delta$ such that $g(N)\cap N\not=
\emptyset$ only for finitely many $g\in \Gamma$. Namely, for
a point $(x,y,z)\in X^3-\Delta$ choose $\nu>0$ sufficiently
small that $\min\{d(x,y),d(x,z),d(z,y)\}\geq 2\nu$.
For this number $\nu$ let $R(\nu)>0$ be as in the
definition of a metrically proper action and let $N$
be the open $e^{-R(\nu)}$-neighborhood of $(x,y,z)$ in $X^3$; then
$N\cap gN\not=\emptyset$ only for finitely many
$g\in \Gamma$ by the second part of our above definition.
Since $X$ does not have isolated points this implies
that the quotient $(X^3-\Delta)/\Gamma$ is a metrizable
Hausdorff space.

For every $g\in \Gamma$,
the fixed point set ${\rm Fix}(g)$ for the action of
$g$ on $X$ is a closed
subset of $X$.  The boundary $A(g)$ of the
open subset $X^3-{\rm Fix}(g)^3$ of $X^3$
is a closed (possibly empty)
nowhere dense subset of $X^3$. By the above observation,
every point $(x,y,z)\in X^3-\Delta$ admits a neighborhood
$N$ which intersects only finitely many of the
sets $A(g)$ $(g\in \Gamma)$. Since $X$ does not
have isolated points, the set
$X^3-\Delta-\cup_{g\in \Gamma}A(g)$
is open and dense in $X^3-\Delta$.
The restriction of the natural projection
\begin{equation}
\pi:T=X^3-\Delta\to Y=(X^3-\Delta)/\Gamma
\end{equation}
to the open dense set $X^3-\Delta-\cup_{g\in \Gamma}A(g)$
is a local homeomorphism.

The involution
$\iota:X^3\to X^3$ defined by $\iota(a,b,c)=(b,a,c)$ is an
isometry with respect to the metric $d$ on $X^3$ induced
from the metric on $X$, and its
fixed point set is contained in the closed
set $\Delta\subset X^3$.
Thus the restriction of $\iota$ to $T$
acts freely, and it commutes
with the diagonal action of $\Gamma$. In particular,
$\iota$ naturally acts on $Y$ as a
continuous involution and the quotient $Z=Y/\iota$ is
a metrizable Hausdorff space.
There is an open dense subset of $T$
such that the restriction of
the natural projection
$\pi_0:T\to Z$ to this set is a local homeomorphism.

For $x\in X$ and
$\epsilon >0$ denote by $B(x,\epsilon)\subset X$ the open ball of
radius $\epsilon$ about $x$. We next  
recall the well known
notion of north-south dynamics for a homeomorphism of $X$.

{\bf Definition:}
A homeomorphism $g$ of $X$
has \emph{north-south dynamics} with respect to an attracting
fixed point $a\in X$ and a repelling fixed point $b\in X-\{a\}$ if
the following is satisfied.
\begin{enumerate}
\item For every $\epsilon
>0$ there is a number $m>0$ such that $g^m
(X-B(b,\epsilon))\subset B(a,\epsilon)$ and $g^{-m}
(X-B(a,\epsilon))\subset B(b,\epsilon)$.
\item There is a number $\delta >0$ such that $\cup_{m\in
\mathbb{Z}}g^m(X-B(a,\delta)-B(b,\delta))=X-\{a,b\}$.
\end{enumerate}
We call $a$ the \emph{attracting} and $b$ the \emph{repelling}
fixed point of $g$, and $(a,b)$ is the \emph{ordered} pair
of fixed points for $g$.

The next definition formalizes the idea
that the dynamics of each element
of a group $G$ of homeomorphisms of $(X,d)$
is uniformly similar to north-south-dynamics on a
metrically large scale.

{\bf Definition:} The action of an arbitrary group $G$
on a metric space $(X,d)$ of finite diameter
is called
\emph{weakly hyperbolic} if for every $\epsilon >0$
there is a number $b=b(\epsilon)\in (0,1)$ with
the following property. Let $x,y\in X$ with $d(x,y)\geq 2\epsilon$
and let $g\in G$ be such that $d(gx,gy)\geq 2\epsilon$.
Let $z\in X-\{x,y\}$ be such
that $\min\{d(gz,gx),d(gz,gy)\}\geq \epsilon$; then $d(gw,gx)\leq
d(z,y)^{b}/b$ for every $w\in X$ with $d(w,x)\leq
\epsilon$.

\bigskip

Let again $\Gamma$ be a countable group which
admits an action on a metric space
$(X,d)$ of finite diameter without isolated points
by homeomorphisms such that
the diagonal action on $T=X^3-\Delta$ is metrically proper.
Using the above notations,
let $C\subset T$ be an open set
whose closure $\overline{C}$ has positive distance
to $\Delta$ and is
mapped by the projection $\pi_0:T\to Z=Y/\iota$ 
homeomorphically into $Z$.
This means that for every $g\in \Gamma$, either $g$ fixes
$\overline{C}\cup \iota \overline{C}$ pointwise or
$g(\overline{C}\cup \iota\overline{C})\cap
(\overline{C}\cup \iota\overline{C})=\emptyset.$
We assume that $C$ is of
the form $C=U\times V\times W$ where $U,V,W\subset X$ are open
and pairwise of positive distance; say that the distance
between any two of these sets is not smaller than a
number $4\nu>0$.
For $R(\nu)>0$
as in the definition of a metrically proper action
we also assume that the diameter of $C$ is smaller
than $e^{-R(\nu)}$. Let ${\cal H}_C$ be the vector
space of all H\"older continuous functions
$f:T\to \mathbb{R}$ supported in $C$.
This means that for every $f\in {\cal H}_C$
there is some $\alpha\in
(0,1)$ and some $q>0$ such that $\vert f(x)-f(y)\vert\leq q
d(x,y)^\alpha$ for all $x,y\in C$.

The following proposition is the analogue
of Lemma 2.1. For its formulation,
denote by ${\cal Q}$ the vector space
of all quasi-morphisms of $\Gamma$.

\bigskip

{\bf Lemma 3.1:} {\it Let $(X,d)$ be a metric
space of finite diameter without isolated
points. Let $\Gamma$ be a countable
group which admits a weakly hyperbolic action
by homeomorphisms of $X$ such that the action of 
$\Gamma$ on $T=X^3-\Delta$ is metrically proper.
Then for every open set $C\subset T$
whose closure 
projects homeomorphically into $Z=(T/\Gamma)/\iota$
there is
a linear map $\Psi:{\cal H}_C\to {\cal Q}$.} 

{\it Proof:} Using the above notations,
write $C=U\times V\times W$ where
$U,V,W\subset X$ are open and pairwise of distance
at least $4\nu>0$. Assume that the diameter of
$C$ is at most $e^{-R(\nu)}$.
The product structure of $T$ defines a natural foliation ${\cal
F}$ on $T$ by requiring that the leaf of ${\cal F}$ through
$(a,b,c)\in T$ equals the set
$F(a,b)=\{(a,b,d)\mid d\in X-\{a,b\}\}$. Thus a leaf
of ${\cal F}$ is determined by two distinct
points in $X$, and the leaf $F(a,b)$ determined by
$a\not= b\in X$ can naturally be identified with
$X-\{a,b\}$. The foliation ${\cal F}$ is
invariant under the action of $\Gamma$ and hence it projects to a
foliation ${\cal F}_0$ on $Y=T/\Gamma$.

Let $\mu_W$ be a
Borel probability measure on $W$ which is positive on
open sets. Choose a nontrivial
H\"older continuous function $\psi:X\times X\to [0,1]$
supported in $U\times V$ and let $\mu_{\cal F}$ be the family
of $\Gamma$-invariant
$\iota$-invariant Borel measures on the leaves of
${\cal F}$ which is determined by the requirement that
for every $(u,v)\in U\times V$ the
restriction of $\mu_{\cal F}$ to $F(u,v)\cap C \sim (u,v)\times W$
equals $\psi(u,v)\mu_W$.

We divide the proof of our lemma into
two steps.
As a notational convention,
for $x\in X$ and $\epsilon >0$ denote
as before by $B(x,\epsilon)$ the
open $\epsilon$-ball about $x$.

{\sl Step 1:}

In our first step we construct for a given H\"older 
continuous function
$f\in {\cal H}_C$ supported in $C$ 
a function $\Psi(f):\Gamma\to
\mathbb{R}$. 
For this recall that by our choice of $C$, 
every $g\in \Gamma$ either fixes $\overline{C}\cup
\iota\overline{C}$ pointwise or we have
$g(\overline{C}\cup \iota
\overline{C})\cap (\overline{C}\cup \iota\overline{C})=\emptyset$.
Therefore every 
function $f\in {\cal H}_C$ 
uniquely determines a continuous
$\Gamma$-invariant $\iota$ anti-invariant function $\tilde f$ on
$T$ which is supported in $\cup_{g\in \Gamma}g(C\cup \iota C)$
and whose restriction to $C$ coincides with $f$.
This means that $\tilde f(\iota x)=-\tilde f(x)$ for all $x\in T$,
$\tilde f(gx)=\tilde f(x)$ for all $g\in \Gamma$ and that moreover
the restriction of $\tilde f$ to $C$ coincides with $f$.
We claim that for all $x\not=y\in X$ and
any neighborhood $A$ of $x$, $B$ of $y$ we have
$\int_{F(x,y)-A-B}\vert \tilde f\vert d\mu_{\cal F}<\infty$
where as before, we identify the leaf $F(x,y)$ of the
foliation ${\cal F}$ with the set $X-\{x,y\}$.

For this
consider first the case that $d(x,y)\geq 2\nu$
where $\nu>0$ is as above determined by the choice of 
the set $C$.
Let $k_0\geq 1$ be the smallest integer which is not
smaller than $-\log \nu$. If $z\in X$ is such that
$d(x,z)\leq e^{-k_0}$ then $d(x,z)=\min\{d(x,z),d(y,z)\}$
and hence by the first requirement in the
definition of a metrically proper action, for every
$k\geq k_0$ the number of
elements $g\in \Gamma$ with the property that
$g(C\cup \iota C) \cap (F(x,y)\cap
(B(x,e^{-k})-B(x,e^{-k-1})))\not=\emptyset$ is bounded
from above by a
constant $m(\nu)>0$ only depending on $\nu$ but not on $k$ and
$(x,y)$. Since $\tilde f$ is invariant under the action of
$\Gamma$ and supported in $\cup_{g\in \Gamma}g(C\cup \iota C)$ and
since the measures $\mu_{\cal F}$ are invariant under the action
of $\Gamma$ this implies that
\begin{equation}
\int_{F(x,y)\cap (B(x,e^{-k})-B(x,e^{-k-1}))}\vert
\tilde f\vert d\mu_{\cal F}\leq
m(\nu)\sup\{\vert f(z)\vert\mid
z\in C\}
\end{equation}
for every $k\geq k_0$. The same estimate also holds for the
integral\\ $\int_{F(x,y)\cap (B(y,e^{-k})-B(y,e^{-k-1}))}\vert \tilde
f\vert d\mu_{\cal F}$ provided that $k\geq k_0$.

Similarly, since the diameter of $X$ is finite the set
$F(x,y)-B(x,e^{-k_0})-B(y,e^{-k_0})$ is the
union of finitely many subsets of the form
\[\{z\mid e^{-k_0+m-1}\leq \min\{d(z,x),d(z,y)\}
\leq e^{-k_0+m}\}\quad (m\geq 1).\] Using once more the
definition of a metrically proper action we conclude that
the number of elements $g\in \Gamma$ such that $g(C\cup
\iota C)\cap (F(x,y)-B(x,e^{-k_0})-B(y,e^{-k_0}))\not=\emptyset$ is
bounded from above by a constant only depending on $\nu$.
In particular, the value of the integral
$\int_{F(x,y)-B(x,e^{-k_0})-B(y,e^{-k_0})}\vert \tilde f\vert
d\mu_{\cal F}$ is bounded from above by a universal multiple
of the supremums norm of $f$. Together
we conclude that
for every neighborhood $A$
of $x$, $B$ of $y$ the integral
$\int_{F(x,y)-A-B}\vert \tilde f\vert d\mu_{\cal F}$
exists, i.e. our claim holds true whenever
$d(x,y)\geq 2\nu$.

Now let $x\not=y\in X$ be
arbitrary points such that the support of $\tilde f$ intersects
the leaf $F(x,y)$. 
Since $\tilde f$ is supported in 
$\cup_{g\in \Gamma} g(C\cup \iota C)$, 
there is then an element $g\in \Gamma$ with
$d(gx,gy)\geq 2\nu$. By invariance of $\tilde f$ and $\mu_{\cal
F}$ under the action of $\Gamma$, for every neighborhood $A$ of
$x$, $B$ of $y$ we have 
\begin{equation}\int_{F(x,y)-A-B}\vert \tilde f\vert
d\mu_{\cal F}= \int_{F(gx,gy)-gA-gB}\vert \tilde f\vert d\mu_{\cal
F}\end{equation}
where $gA,gB$ is a neighborhood of $gx,gy$. Thus
indeed
$\int_{F(x,y)-A-B}\vert f\vert d\mu_{\cal F}
<\infty$ for any two points
$x\not= y\in  X$ and any neighborhoods
$A$ of $x$, $B$ of $y$ which shows our above claim.

Recall that $C=U\times V\times W$ for open disjoint subsets
$U,V,W$ of $X$. Fix a point $x\in U$ and let $A\subset U$ be a
small closed metric ball centered at $x$. For 
$f\in {\cal H}_C$ and $g\in \Gamma$ such that
$gx\not= x$ define
\begin{equation}\label{integral}
\Psi(f)(g)=\int_{F(x,gx)-A-g(A)}\tilde fd\mu_{\cal F}
\end{equation}
and if $gx=x$ then define
$\Psi(f)(g)=0$.
By our above consideration, the integral (\ref{integral}) exists
and hence it defines a function
$\Psi(f):\Gamma\to \mathbb{R}$.
Moreover, the assignment $\Psi:f\in {\cal H}_C\to \Psi(f)$ is
a linear map from the vector space ${\cal H}_C$ into the
vector space of all functions on $\Gamma$.

{\sl Step 2:}

In a second step, we
show that for every $f\in {\cal H}_C$ the function
$\Psi(f):\Gamma\to \mathbb{R}$ is a quasi-morphism,
i.e. that we have
$\sup_{g,h}\{\vert \Psi(f)(g)+\Psi(f)(h)-\Psi(f)(gh)\vert\} <\infty$.
Observe that by invariance
under $\Gamma$, for $g,h\in
\Gamma$ we have
\begin{align}
\Psi(f)(g)+\Psi(f)(h)-\Psi(f)(gh)= &
\int_{F(x,gx)-A-gA}\tilde f d\mu_{\cal F} \\ +
\int_{F(gx,ghx)-gA-ghA}\tilde fd\mu_{\cal F}- &
\int_{F(x,ghx)-ghA-A} \tilde fd\mu_{\cal F}.\notag
\end{align}
Since $f$ is anti-invariant under
the involution $\iota$ it is
therefore enough to show that there is a number
$c(\nu,f)$ only
depending on $\nu$ and the H\"older norm of $f$ with the
following property.
Let $(x_1,x_2,x_3)\in T$ and let $A_i$ be any neighborhood of
$x_i$ in $X$ $(i=1,2,3)$; then
\begin{equation}
\vert \int_{F(x_1,x_2)-A_1-A_2}\tilde f d\mu_{\cal F} +
\int_{F(x_2,x_3)-A_2-A_3}\tilde fd\mu_{\cal F} +
\int_{F(x_3,x_1)-A_3-A_1} \tilde fd\mu_{\cal F}\vert <c(\nu,f).
\end{equation}

For this recall that for $g,h\in \Gamma$
the sets $g\overline{C} , h\overline{C},g(\overline{\iota C}),
h(\overline{\iota C})$
either coincide
or are disjoint. Moreover, if $\tilde f\vert F(y,z)\not\equiv
0$ for some $y\not=z\in X$
then there is some $g\in \Gamma$ such that $d(gy,gz)\geq 2\nu$.
Define ${\cal G}=\{g\in \Gamma\mid \max_{i,j\leq 3}d(gx_i,gx_j)\geq
2\nu\}$ and
\begin{equation} {\cal G}_0=\{g\in {\cal G}\mid
\min\{d(gx_i,gx_j)\mid i\not=j\in \{1,2,3\}
\}\geq \nu\}\end{equation}
and for $i=1,2,3$ define
\begin{equation}
{\cal G}_i=\{g\in {\cal G}\mid
d(gx_i,gx_{i+1})<\nu\}\end{equation}
(indices are taken modulo 3). By the triangle inequality
and the definition of the set ${\cal G}$, the
sets ${\cal G}_i (i=0,\dots,3)$ are pairwise disjoint and their
union equals ${\cal G}$.

If $g\in \Gamma$ is such that $(C\cup \iota C)\cap
gF(x_i,x_{i+1})\not= \emptyset$ then $d(gx_i,gx_{i+1})\geq 2\nu$ and
therefore $g\in {\cal G}_0$ if $\min\{d(gx_{i-1},gx_i),
d(gx_{i-1},gx_{i+1})\}
\geq \nu,\, g\in {\cal G}_{i-1}$ if $d(gx_{i-1},gx_{i})<\nu$ and
$g\in {\cal G}_{i+1}$ otherwise (where indices are
again taken modulo 3). Thus by invariance of 
$\tilde f$ and $\mu_{\cal F}$
under the action of $\Gamma$ and by
the fact that an element $g\in \Gamma$
either fixes $C$ pointwise or is such that
$g C\cap C=\emptyset$ it is enough to show that
there is a number $c_1(\nu,f)>0$
only depending on $\nu$ and the H\"older norm of $f$ such that
for $i=0,\dots,3$ we have
\begin{equation}\label{sum}
\sum_{g\in {\cal G}_{i}}\vert \sum_{j=1}^3
\int_{g(F(x_j,x_{j+1})-A_j-A_{j+1})
\cap (C\cup \iota C)}
\tilde f d\mu_{\cal F}\vert \leq c_1(\nu,f).
\end{equation}

We first establish the estimate (\ref{sum}) for $i=0$.
The case ${\cal G}_0=\emptyset$ is trivial,
so assume that there is some 
$h\in {\cal G}_0$ with the
additional property that
$h(F(x_j,x_{j+1})-A_j-A_{j+1})\cap C
\not=\emptyset$ for some $j\in \{1,2,3\}$.
Recall that $C=U\times V\times W$ where
the diameter of the sets $U,V,W$ is at most
$e^{-R(\nu)}<\nu/4$. Let $Z\subset X$ be the set
of all points whose distance to $U\cup V$ is
at least $\nu$. Then $h(x_j,x_{j+1},x_{j+2})\in
U\times V\times Z$ and therefore if $u\in {\cal G}_0$
is another element with 
$u(F(x_j,x_{j+1})-A_j-A_{j+1})
\cap C\not=\emptyset$ 
then $uh^{-1}(U\times V\times Z)
\cap U\times V\times Z\not=\emptyset$.
Using
the second property in the definition of a metrically proper
action we conclude that
the number of elements $u\in {\cal G}_0 $
with this property
is bounded from above by a constant
only depending on $\nu$. The same argument
also applies to elements $g\in {\cal G}_0$ with
$g(F(x_j,x_{j+1})-A_j-A_{j+1})\cap \iota C\not=
\emptyset$ for some $j\in\{1,2,3\}$. 
As a consequence, for $i=0$ the number of 
nonzero terms in the sum (\ref{sum}) is
bounded from above by a universal constant and
the estimate (\ref{sum}) holds true for $i=0$.
Thus by symmetry in $i\in \{1,2,3\}$ and by
invariance under the action of $\Gamma$
it now suffices to show the estimate (\ref{sum})
for $i=3$.

By definition, for $g\in {\cal
G}_3$ we have $d(gx_1,gx_3)<\nu$ and therefore $gF(x_1,x_3)\cap
(C\cup \iota C) =\emptyset$. This means that
\begin{align}\label{doublesum}
\sum_{g\in {\cal G}_3}\vert \sum_{j=1}^3
\int_{g(F(x_j,x_{j+1})-A_j-A_{j+1})
\cap (C\cup \iota C)} & \tilde fd\mu_{\cal F}\vert \\
=\sum_{g\in {\cal G}_3}\vert
\int_{g(F(x_1,x_2)-A_2-A_3)\cap (C\cup \iota C)}\tilde fd\mu_{\cal F}
+ & \int_{g(F(x_2,x_3)-A_3-A_1)\cap (C\cup \iota C)}
\tilde fd\mu_{\cal F}\vert. \notag
\end{align}

By assumption, the action of
$\Gamma$ on $X$ is weakly hyperbolic. Thus there is a constant
$b>0$ depending on $\nu$
such that for all $(x,y)\in X\times X$ with
$d(x,y)\geq 2\nu$,
for all $k\geq -\log \nu$ and for all
$z\in X$ with $d(z,x)\leq \nu$
we have $d(gz,gx)\leq e^{-kb}/b$
whenever $g\in \Gamma$ is such that
$(C\cup \iota C)\cap g(F(x,y)\cap B(y,e^{-k}))\not=\emptyset$.
In particular, for every $w\in W$ the
distance between
$(gx,gy,w)$ and $(gz,gy,w)=\iota(gy,gz,w)$
is at most $e^{-kb}/b$.

Now $\tilde f$ is a $\Gamma$-invariant $\iota$-anti-invariant
function on $T$ which is supported in
$\cup_{g\in \Gamma}g(C\cup \iota C)$ and
whose restriction $f$ to $C$ satisfies
$\vert f(v)-f(w)\vert \leq qd(v,w)^\alpha$ 
for some $\alpha >0,q>0$ and
for all $v,w\in C$. Moreover,
$\mu_{\cal F}$ is a
$\iota$-invariant $\Gamma$-invariant family
of measures on the leaves of ${\cal F}$ whose restriction to $C$
is of the form $\psi\mu_W$ for a H\"older continuous
function $\psi$ supported in  $U\times V$.
As a consequence of our above discussion on the
effect of weak hyperbolicity
we conclude that there
is a number $\beta >0$ depending only
on $\nu$ and the H\"older norm of $f$
with the
following property.
Let $x,y\in X$ with $d(x,y)\geq 2\nu$;
if for some $k\geq -\log \nu$ the element
$g\in \Gamma$ is such that $(C\cup \iota C)\cap
g(F(x,y)\cap (B(y,e^{-k})-B(y,e^{-k-1})))\not=\emptyset$ then
for every $z\in X$ with $d(z,x)\leq \nu$
and every neighborhood $A$ of $y$ we have
\begin{equation}\label{intsum}\vert \int_{g(F(x,y)-A)\cap
(C\cup \iota C)}\tilde fd\mu_{\cal F}+
\int_{g(F(y,z)-A)\cap (C\cup \iota C)}\tilde f
d\mu_{\cal F}\vert \leq
e^{-k\beta}/\beta.\end{equation}

By the first property in the
definition of a metrically proper action there is a constant
$c>0$ only depending on $\nu$ such that for all $x,y\in X$
with $d(x,y)\geq 2\nu$
and every $k\geq -\log \nu$
there are at most $c$ elements $g\in \Gamma$
with $(C\cup \iota C)\cap g(F(x,y)\cap (B(y,e^{-k})-B(y,e^{-k-1})))\not=
\emptyset$.
Together with the estimate (\ref{intsum}) we conclude that
there is a constant $c_2(\nu,f)>0$ which only
depends on $\nu$ and on the H\"older norm of $f$ with
the following property.
For $x,y\in X$ with $d(x,y)\geq 2\nu$,
every $z\in X-\{x,y\}$ with $d(x,z)\leq \nu$
and every neighborhood $A$ of $y$
we have
\begin{equation}\label{extended} 
\sum_{\{g\in \Gamma\mid d(gx,gz)<\nu\}}
\vert \int_{g((F(x,y)-A)\cup (F(y,z)-A))\cap (C\cup \iota C)} 
\tilde fd\mu_{\cal F}\vert <c_2(\nu,f).
\end{equation}
Now if $g\in {\cal G}_3$ is such that
$gF(x_1,x_2)\cap (C\cup \iota C)\not=\emptyset$ then 
with $y_i=gx_i$ we have
$d(y_1,y_3)<\nu, d(y_1,y_2)\geq 2\nu$. 
For every other $h\in {\cal G}_3$ with
$h F(x_1,x_2)\cap (C\cup \iota C)\not=\emptyset$ 
we obtain $d(hg^{-1}y_1,hg^{-1}y_2)\geq 2\nu$ and 
$d(hg^{-1}y_1,hg^{-1}y_3)\leq \nu$. 
By invariance of $\tilde f$ and
$\mu_{\cal F}$ under the action
of $\Gamma$, inequality (\ref{sum}) above now follows
from the estimate (\ref{extended}) and the
equation (\ref{doublesum}).

As a consequence, for every $f\in {\cal H}_C$ the function
$\Psi(f)$ on $\Gamma$ is indeed a quasi-morphism. 
By construction, the assignment $f\to \Psi(f)$ is 
moreover linear. This completes the proof of the
lemma.
\qed

\bigskip

In Lemma 3.1 we constructed a linear
map $\Psi$ from the vector space ${\cal H}_C$ onto
a vector space $\Psi({\cal H}_C)\subset {\cal Q}$ 
of quasi-morphisms for the group $\Gamma$.
It follows from our construction that
for a suitable choice of our set $C$ the vector space
$\Psi({\cal H}_C)$ is infinite dimensional. 
As in Section 2, the map $\Psi$ then induces 
via composition with the natural projection
a linear map
$\Theta:{\cal H}_C\to H_b^2(\Gamma,\mathbb{R})$.
However, a priori the image of $\Theta$ 
may be trivial or finite dimensional. To establish 
that the subspace of $H_b^2(\Gamma,\mathbb{R})$ obtained
in this way as the sets $C$ vary,  
we use an additional assumption on $\Gamma$ 
which is motivated by the work \cite{BF} of
Bestvina and Fujiwara. For this recall that
a homeomorphism $g$ of $X$ which acts with 
north-south dynamics has an attracting fixed
point $a\in X$ and a repelling fixed point $b\in X-\{a\}$.
We call $(a,b)$ the \emph{ordered} pair of fixed
points for $g$. 
We show.

\bigskip

{\bf Proposition 3.2:} {\it In the situation 
described in Lemma 3.1, assume in addition
that the group $\Gamma$ contains a free subgroup 
$G$ with two generators and the following properties.
\begin{enumerate} \item Every $e\not=g\in G$ acts with
north-south dynamics on $X$.
\item There are infinitely many $g_i\in G$ $(i>0)$
such that the $\Gamma$-orbits of the ordered
pairs of fixed points of the elements $g_i,g_j^{-1}$ $(i,j>0)$ 
are pairwise disjoint.
\end{enumerate}
Then the images of the spaces ${\cal H}_C$ under the map $\Theta$
for suitable choices of $C\subset T$ 
span an infinite dimensional subspace of 
$H_b^2(\Gamma,\mathbb{R})$.}

{\it Proof:} Continue to use the assumptions
and notations from Lemma 3.1 and its proof.
We have to show that the
bounded cohomology classes $\Theta(f)$ $(f\in {\cal H}_C)$
defined
by the quasi-morphisms $\Psi(f)$ constructed in Lemma 3.1
for suitable choices of $C$ span an infinite dimensional
subspace of the kernel
of the map $H_b^2(\Gamma,\mathbb{R})\to H^2(\Gamma,\mathbb{R})$.
For this let $G$ be the free subgroup of $\Gamma$ with
two generators as in the statement of our proposition.

Let $g,h\in G-\{e\}$ be such that the 
$\Gamma$-orbit of the ordered
pair $(a,b)$ of fixed points for $g$ is distinct from the
$\Gamma$-orbit of the ordered pair $(a^\prime,b^\prime)$ of 
fixed points for $h$. Then the leaves $F(a,b),
F(a^\prime,b^\prime)$
of the foliation ${\cal F}$ project to distinct leaves
$L,L^\prime$ of the induced foliation ${\cal F}_0$ on 
$Y=T/\Gamma$. We
claim that the closures of these leaves do not intersect.
For this
denote as before by $\pi:T\to Y$ the natural projection.
Let $\epsilon_0 >0$ be
sufficiently small that $d(\{a,b\},\{a^\prime,b^\prime\})\geq
2\epsilon_0$. Since
$g,h$ act on $X$ with north-south dynamics and fixed points $a,b$
and $a^\prime,b^\prime$, there is a number $\epsilon <\epsilon_0$
with the property that the projection $\pi$ maps the set
$\{(a,b,x)\mid d(x,\{a,b\})\geq \epsilon\}$ onto $L$ and that
$\pi$ maps $\{(a^\prime,b^\prime,y)\mid d(y,\{a^\prime,b^\prime\})
\geq \epsilon\}$ onto
$L^\prime$.

Assume to the contrary that the closures of the leaves
$L,L^\prime$ contain a common point. By our above observation,
this implies that there is a sequence $(g_i)\subset \Gamma$ of
pairwise distinct elements and there are sequences $(x_i)\subset X,
(y_i)\subset X$ such that
\[d(x_i,\{a,b\})\geq \epsilon, d(y_i,\{a^\prime,b^\prime\})
\geq \epsilon\quad \hbox{\rm for all }\, i\] and that
$d(g_i(a,b,x_i),(a^\prime, b^\prime,y_i))\to 0$. In particular,
for every $\delta >0$ there are infinitely many distinct elements
$u\in \Gamma$ such that $d(a^\prime,ua)< \delta,d(b^\prime,ub)<
\delta$ and that $u(X-B(a,\epsilon)-B(b,\epsilon))\cap
X-B(a^\prime,\epsilon)-B(b^\prime,\epsilon)\not=\emptyset$.
However this contradicts the
second requirement in the definition of a
metrically proper action.
As a consequence, the closures of the leaves $L,L^\prime$
in $Y$ are disjoint.

Let $g\in G$ and let $a$ be the attracting and $b$ be the
repelling fixed point of $g$. 
Choose the set $C=U\times V\times W\subset T$
as in Lemma 3.1 and its proof in such a way that
$a\in U,b\in V$. This is possible since 
the action of $\Gamma$ on $X$ is metrically proper and
hence the stabilizer of $\{a,b\}$ in $\Gamma$ 
acts freely on an open subset of 
$X-\{a,b\}$. Let $x\in U-\{a\}$ and choose a closed
neighborhood $A\subset U-\{a\}$ of $x$ for
the construction of the quasi-morphism $\Psi(f)$.
Since $g$ acts on $X$ with north-south
dynamics there is a closed subset $D$ of $X-\{a,b\}$ 
with dense interior whose
distance to $\{a,b\}$ is positive and which is a fundamental
domain for the action on  $X-\{a,b\}$ of the cyclic subgroup of
$G$ generated by $g$. For the measures $\mu_{\cal F}$ on the
leaves of the foliation ${\cal F}$ as 
in the proof of Lemma 3.1 we may assume that the
support of $\mu_{\cal F}$ intersects
$F(a,b)$ and that the
$\mu_{\cal F}$-mass of the boundary of $D$ viewed as a subset of
$F(a,b)$ vanishes. 
Let $f\in {\cal H}_C$ and let $\tilde f$ be
the $\Gamma$-invariant $\iota$-anti-invariant
function on $T$ defined by $f$ as in the proof of
Lemma 3.1.
By the discussion in Step 1 of the proof of
Lemma 3.1, the integral
$\int_D\tilde fd\mu_{\cal F}$ exists. Let $\Psi(f)$
be the quasi-morphism of $\Gamma$ defined by $f$ as
in (8) of Lemma 3.1. We claim that
\begin{equation}
\lim_{k\to \infty}\Psi(f)(g^k)/k=\int_D\tilde fd\mu_{\cal F}.
\end{equation}

To show our claim,
observe that as $k\to \infty$ the
diameter of the sets $g^kA$ tends to $0$ and $g^kx\to a\in X-A$.
Choose a small closed ball $B\subset V$ about $b$. By the
consideration in the proof of Lemma 3.1,
for sufficiently large $k$ the absolute value of
the difference
\begin{equation}
\int_{F(g^kx,b)-g^kA-B}\tilde f d\mu_{\cal F}-
\int_{F(g^kx,x)-g^kA-A}\tilde fd\mu_{\cal F}
\end{equation}
is bounded from above by a
constant not depending on $k$.
As a consequence, it is enough to
show that
\begin{equation}
\int_{F(g^kx,b)-g^kA-B}\tilde fd\mu_{\cal F}/k
\to \int_D \tilde f d\mu_{\cal F} \quad (k\to\infty),
\end{equation}
and this in turn is equivalent to
\begin{equation}
\int_{F(x,b)-A-g^{-k}B}\tilde fd\mu_{\cal F}/k
\to \int_D \tilde f d\mu_{\cal F}\quad (k\to\infty).
\end{equation}

Choose in particular
$B=\{b\}\cup_{j\leq 0}g^jD$. Then $B-g^{-k}B=\cup_{j =0}^{k-1}
g^{-j} D$ for every $k>0$. Thus
for every small
ball $E\subset X-B$
about the attracting fixed point $a$ for $g$ we have
\begin{equation}
\lim_{k\to \infty}
\int_{F(x,b)-A-g^{-k}B}\tilde f
d\mu_{\cal F}/k=\lim_{k\to\infty}
\int_{F(a,b)-E-g^{-k}B}\tilde fd\mu_{\cal F}/k
=\int_D\tilde fd\mu_{\cal F}.
\end{equation}
This shows our above claim.

Let again $g\in G$ with attracting fixed point $a\in X$,
repelling fixed point $b\in X-\{a\}$ and assume that
the ordered pair $(a,b)$ is not contained in 
the $\Gamma$-orbit of the ordered pair $\iota(a,b)=(b,a)$.
By our above consideration, the closure
of the projection of the leaf $F(a,b)$ to $Y$
is disjoint from the closure of the projection 
of $\iota F(a,b)=F(b,a)$. As before, let 
$D\subset F(a,b)$ be a closed fundamental domain
for the action on $X-\{a,b\}\sim F(a,b)$ 
of the cyclic group generated by $g$. By the
second requirement in the definition of a metrically
proper action, there are
only finitely many $h\in \Gamma$ with $hD\cap D\not=\emptyset$.
Denote by $\pi:T\to Y$ the canonical projection.
The measures
$\mu_{\cal F}$ project to a family of measures on the leaves of
the foliation ${\cal F}_0=\pi{\cal F}$ on $Y$.
For $f\in {\cal H}_C$ the function
$\tilde f$ projects to a function $f_0$ on $Y$. 
Since $hD\cap D\not=\emptyset$ for only finitely
many $h\in \Gamma$, 
the integral $\int_D\tilde fd\mu_{\cal F}$ is a positive
bounded multiple of the integral $\int_{\pi F(a,b)}f_0d\mu_0$.
By our above consideration,
the closure $L$ of the projection of the leaf
$F(a,b)$ to $Y$ is disjoint from the closure of its
image $F(b,a)$ under the involution $\iota$ and therefore
for any given number $q\in \mathbb{R}$ there is
a H\"older function $f\in {\cal H}_C$ such that the quasi-morphism
$\Psi(f)$ defined as above by $f$ satisfies
$\lim_{k\to \infty} \Psi(f)(g^k)/k=q$.

By our assumption, there are infinitely
many elements $g_i\in G$ $(i>0)$ which act on $X$ with north-south
dynamics and such that the ordered pairs of 
fixed points of $g_i,g_j^{-1}$ are pairwise
contained in distinct $\Gamma$-orbit on $X$.
In particular,
for $i\not=j$ the closures of the
projections to $Y$ of the leaves of the foliation ${\cal F}$ which
are determined by the fixed points of $g_i,g_j$ are disjoint.
Now for any finite set $\{h_1,\dots, h_m\}\subset \{g_i\mid
i>0\}\subset G$ choose the set $C$ as
above in such a way that it intersects each of the leaves of
${\cal F}$ determined by the 
ordered pair of fixed points of $h_i$; this can
always be achieved by allowing for our construction a set $C$
which consists of finitely many components satisfying each our
above assumptions. Our above discussion shows that for an
arbitrarily chosen collection
$\{q_1,\dots,q_m\}\subset \mathbb{R}$ of real numbers there is
a suitable choice of the
function $f\in {\cal H}_C$ so that the
quasi-morphism $\Psi(f)$ for $\Gamma$ defined by $f$ satisfies
$\lim_{k\to \infty}\Psi(f)(h_i^k)/k=q_i$ for $1\leq i\leq k$.

For $f\in {\cal H}_C$ the cohomology class
$\Theta(f)\in H_b^2(\Gamma,\mathbb{R})$
vanishes if and only if there is a
homomorphism $\eta\in H^1(\Gamma,\mathbb{R})$
such that $\sup_{g\in \Gamma}\vert \Psi(f)(g)-\eta(g)\vert <\infty$
(compare the discussion in Section 2). This homomorphism
then restricts to a homomorphism of the group $G$.
Now $G$ is a free group with two generators and hence we have
$H^1(G,\mathbb{R})=\mathbb{R}^2$. More precisely,
if $u_1,u_2$ are such free generators for $G$ then
every homomorphism $\eta:G\to \mathbb{R}$ is
determined by its value on $u_1,u_2$.
In particular,
for any finite subset $\{h_1,\dots,h_m\}\subset G$
there are two elements in this collection, say the
elements $h_1,h_2$, such that for
every quasi-morphism $\eta$ for $G$ which is equivalent
to a homomorphism and every $j\in \{3,\dots,m\}$
the quantity
$\lim_{k\to\infty}\eta(h_j^k)/k$ is uniquely determined by
$\lim_{k\to\infty}\eta(h_i^k)/k$ $(i=1,2)$.
Together with the above observation that for any finite
subset $\{h_1,\dots,h_m\}$ of $\{g_i\mid i> 0\}$
we can find a quasi-morphism for $\Gamma$ for which these
limits assume arbitrarily prescribed values
we conclude that there are infinitely
many quasi-morphisms for $\Gamma$ whose restrictions to $G$ define
linearly independent elements of $H_b^2(G,\mathbb{R})$.
This shows that
the kernel of the map $H_b^2(\Gamma,\mathbb{R})\to
H^2(\Gamma,\mathbb{R})$ is infinite dimensional and
completes the proof of the proposition.
\qed

\bigskip

{\bf Remark:} Our above proof also shows the following.
Let $\Gamma$ be a countable group which admits
a weakly hyperbolic action by homeomorphisms
of a metric space $X$ of finite diameter such that
the action of $\Gamma$ on $T=X^3-\Delta$ is metrically
proper. Let $g_i\in \Gamma$  
be elements which act
with north-south dynamics on $X$ with ordered pairs of fixed
points $(a_i,b_i)$ $(i=1,\dots,k)$. If the
$\Gamma$-orbits of $(a_i,b_i),(b_j,a_j)$ $(i,j\leq k)$
are all disjoint then for every tuple $(q_1,\dots,q_k)\in
\mathbb{R}^k$ there is a quasi-morphism
$\phi$ for $\Gamma$ with $\lim_{\ell\to \infty}
\phi(g_i^\ell)/\ell=q_i$ for every $i\leq k$.  

\bigskip

The following theorem 
is the main technical result of 
this note. For its formulation, recall that 
the free group $G$ with two generators is the
fundamental group of a convex cocompact hyperbolic surface whose
limit set $B$ is just the \emph{Gromov boundary} of $G$.

\bigskip

{\bf Theorem 3.3:} {\it Let $(X,d)$ be a
metric space of finite diameter without isolated
points. Let $\Gamma$ be a countable
group which admits a weakly hyperbolic action by
homeomorphisms of $X$.
Assume that $\Gamma$ contains a free subgroup
$G$ with two generators and the following properties.
\begin{enumerate} \item Every $e\not=g\in G$ acts with
north-south dynamics on $X$.
\item There are infinitely many $g_i\in G$ $(i>0)$
such that the $\Gamma$-orbits of the ordered pairs
of fixed points of the elements $g_i,g_j^{-1}$ $(i,j>0)$ 
are pairwise disjoint.
\item There is a $G$-equivariant continuous embedding
of the Gromov boundary of $G$ into $X$.
\end{enumerate}
If the action of $\Gamma$ on the space of triples of pairwise
distinct points in $X$ is metrically proper
then for every $p\in (1,\infty)$ 
the kernel of the map
$H_b^2(\Gamma,\ell^p(\Gamma))\to H^2(\Gamma,\ell^p(\Gamma))$ is
infinite dimensional.}

{\it Proof:} Let $\Gamma$ be a countable group acting by
homeomorphisms on a
metric space $(X,d)$ of finite diameter without isolated
points.
Assume that the action of $\Gamma$ is weakly hyperbolic and that
the diagonal action of $\Gamma$
on the space $T=X^3-\Delta$ of
triples of pairwise distinct points in $X$ is metrically proper.
Write $Y=T/\Gamma$ and 
denote as before by $\iota:T\to T$ the natural involution
which exchanges the first two points in a triple.
Let $G$ be a free subgroup of $\Gamma$ with two
generators as in the statement of the theorem. In particular,
we assume that there is a continuous
$G$-equivariant embedding of the Gromov
boundary $B$ of $G$ into $X$. 
We have to show that for
every $p\in (1,\infty)$ the
kernel of the map $H_b^2(\Gamma,\ell^p(\Gamma))\to
H^2(\Gamma,\ell^p(\Gamma))$ is infinite dimensional.

Denote by $\Vert \,\Vert_p$ the norm of the Banach space
$\ell^p(\Gamma)$. We assume that $\Gamma$ acts on
$\ell^p(\Gamma)$ by right
translation, i.e. for every $g\in \Gamma$ and every function
$\psi\in \ell^p(\Gamma)$ we have $(g\psi)(h)=\psi(hg)$.
Define an \emph{$\ell^p(\Gamma)$-valued quasi-morphism} for
$\Gamma$ to be a map $\eta:\Gamma\to \ell^p(\Gamma)$ such that
\begin{equation}
\sup_{g,h\in \Gamma}\Vert \eta(g)+g\eta(h)
-\eta(gh)\Vert_p<\infty.
\end{equation}
Two such quasi-morphisms $\eta,\eta^\prime$ are called
equivalent if $\eta-\eta^\prime$ is bounded as a function from
$\Gamma$ to $\ell^p(\Gamma)$, i.e. if there is
a number $c>0$ such that $\Vert (\eta-\eta^\prime)(g)\Vert_p\leq c$
for all $g\in \Gamma$.

By Corollary 7.4.7 in \cite{M}, the cohomology group
$H_b^2(\Gamma,\ell^p(\Gamma))$ coincides with the second
cohomology group of the complex
\begin{equation}
0\to L^\infty(\Gamma,\ell^p(\Gamma))^\Gamma\xrightarrow{d}
L^\infty(\Gamma^2,\ell^p(\Gamma))^\Gamma
\xrightarrow{d} L^\infty(\Gamma^3,\ell^p(\Gamma))^\Gamma\to \dots
\end{equation}
with the usual homogeneous coboundary operator $d$.
Let $\psi:\Gamma^2\to \ell^p(\Gamma)$
be any (unbounded) $\Gamma$-equivariant
map; this means that $\psi(hg_1,hg_2)=h(\psi(g_1,g_2))$ for all
$g_1,g_2,h\in \Gamma$.
If the image $d\psi$ of $\psi$ under the
coboundary operator $d$ is \emph{bounded}, then as in the case
of real coefficients, the map $\psi$ defines a class
in the kernel of the natural map
$H_b^2(\Gamma,\ell^p(\Gamma))\to H^2(\Gamma,\ell^p(\Gamma))$.
Let $e$ be the unit element in
$\Gamma$ and
define a map $\phi:\Gamma\to \ell^p(\Gamma)$
by $\phi(v)=\psi(e,v)$.
Then for $g,h,u\in \Gamma$ we have
$d\psi(g,h,u)=\psi(h,u)-\psi(g,u)+\psi(g,h)=h\phi(h^{-1}u)-
g\phi(g^{-1}u)+g\phi(g^{-1}h)=g(\phi(g^{-1}h)+g^{-1}h\phi(h^{-1}u)-
\phi(g^{-1}u))$. Since $\Gamma$ acts isometrically on
$\ell^p(\Gamma)$, we conclude that
$d\psi$ is bounded if and only if $\phi$ defines an
$\ell^p(\Gamma)$-valued quasi-morphism for $\Gamma$.
Now by equivariance, $\psi$ is uniquely determined
by $\phi$ and therefore
every equivalence class of an
$\ell^p(\Gamma)$-valued quasi-morphism gives rise to a cohomology
class in the kernel of the natural map
$H_b^2(\Gamma,\ell^p(\Gamma))\to H^2(\Gamma,\ell^p(\Gamma))$.
This cohomology class vanishes if and only if
there is a map $\eta:\Gamma\to \ell^p(\Gamma)$
which satisfies $\eta(gh)=\eta(g)+g\eta(h)$ for all
$g,h\in \Gamma$ and
such that $\phi-\eta$ is bounded.

Let again $T$ be the space of triples of pairwise
distinct points in $X$. The group $\Gamma$ and the involution
$\iota$ act on $T$, and these
actions commute; we denote as before by $Z$ the corresponding
quotient.
As above, let $C\subset T$ be
a set of positive distance to $\Delta$ and sufficiently small
diameter which is mapped homeomorphically into
the quotient $Z$.

Let $\hat T=T\times
\Gamma$
and define $\hat{\cal H}$ to be the vector space of
all functions $f:\hat{T}\to \mathbb{R}$ supported in
$C\times \Gamma$ with the following property.
For $g\in \Gamma$ write $f_g(x)=f(x,g)$; we view
$f_g$ as a function $C\to \mathbb{R}$. We require
that there is some $\alpha\in
(0,1)$ such that the H\"older-$\alpha$-norms $\Vert
f_g\Vert_\alpha$ of the functions $f_g$ $(g\in \Gamma)$
on $C$ satisfy
$\sum_{g\in \Gamma}\Vert f_g\Vert_\alpha^p<\infty$.
Then for each $y\in C$ the function $f_y:g\to f(y,g)$ is contained
in $\ell^p(\Gamma)$ and therefore
the assignment $y\in C\to f_y$ defines a
(H\"older continuous) map of $C$ into
$\ell^p(\Gamma)$. The set of all such functions naturally has the
structure of an infinite dimensional vector space.
Extend the function $f\in \hat{\cal H}$ to a
function $\hat f$ on $\hat T$ which is
anti-invariant
under the involution $\iota:(\zeta,g)=(\iota\zeta,g)$ and
satisfies $\hat f(gz,u)=f(z,ug)$ for all
$z\in T$, all $g,u\in \Gamma$.

As in the proof of Lemma 3.1 above, 
assume that $C=U\times V\times W$
for open subsets $U,V,W$ of positive distance and
sufficiently small diameter. 
Recall from the proof of Lemma 3.1 the
definition of the foliation ${\cal F}$ of 
$T$ and the measures $\mu_{\cal F}$.
Choose a small closed ball
$A\subset U$, a point $x\in A$ and for
$g\in \Gamma$ define a function
$\Psi(f)(g):\Gamma\to\mathbb{R}$ by
\begin{equation}
\Psi(f)(g)(u)=\int_{F(x,gx)-A-gA}\hat f(y,u)d\mu_{\cal F}(y).
\end{equation}
It follows from our choice of $f$ and the consideration
in Step 1 of the proof of Lemma 3.1
that $\Psi(f)(g)\in \ell^p(\Gamma)$.
On the other hand, by the definition of the
function $\hat f$ we have
\begin{equation}
\int_{F(gx,ghx)-gA-ghA}\hat f(y,u)d\mu_{\cal F}(y)=
\int_{F(x,hx)-A-hA}\hat f(y,ug)d\mu_{\cal F}(y)=\Psi(f)(h)(ug)
\end{equation}
and consequently the estimates in Step 2 
of the proof of Lemma 3.1 show
that the map $\Psi(f)$ is an $\ell^p(\Gamma)$-valued
quasi-morphism for $\Gamma$.
In other words, as in the case of real coefficients we
obtain a linear map $\Theta$ from the vector space
$\hat{\cal H}$ into
the kernel of the natural map
$H_b^2(\Gamma,\ell^p(\Gamma))\to H^2(\Gamma,\ell^p(\Gamma))$
which assigns to a function $f\in \hat{\cal H}$
the cohomology class of the $\ell^p(\Gamma)$-valued
quasi-morphism $\Psi(f)$.

Our goal is to show that the image of the map
$\Theta$ is infinite dimensional. For this
let $G<\Gamma$ be the free
group with two generators as in the statement of the
theorem. Then every function $u\in \ell^p(\Gamma)$ restricts
to a function $Ru\in \ell^p(G)$, and for $g\in G$ we have
$R(gu)=g(Ru)$. Thus for every $f\in \hat{\cal H}$ the
map $\Psi(f):\Gamma\to \ell^p(\Gamma)$ restricts to
an $\ell^p(G)$-valued
quasi-morphism $R\Psi(f):G\to \ell^p(G)$ which defines
a cohomology class $R\Theta(f)\in H_b^2(G,\ell^p(G))$.
If the cohomology class $\Theta(f)$ vanishes
then the same is true for the cohomology
class $R\Theta(f)$.
Thus it is enough to show
that the subspace $\{R\Theta(f)\mid f\in \hat{\cal H}\}$
of $H_b^2(G,\ell^p(G))$ is infinite dimensional.

For this let $B$ be the Gromov boundary of the
free group $G$; this boundary is a Cantor
set on which the group $G$ acts as a group of
homeomorphisms with north-south dynamics.
Assume that
there is a $G$-equivariant continuous
embedding $\rho_0:B\to X$. If we denote by $BT$ the
space of triples of pairwise distinct points in $B$
then the map $\rho_0$ induces a continuous
$G$-equivariant embedding $\rho:BT\to T$.
In the sequel we identify $BT$ with its image
under $\rho$, i.e. we suppress the map
$\rho$ in our notations.
Let $f\in \hat{\cal H}$; for a triple $(x_1,x_2,x_3)\in BT$
and $u\in G$ define
\begin{align}\nu(f)(x_1,x_2,x_3)(u)= &
\int_{F(x_1,x_2)}\hat f(y,u)d\mu_{\cal F}(y)\\+
\int_{F(x_2,x_3)}\hat f(y,u)d\mu_{\cal F}(y) &
+\int_{F(x_3,x_1)}\hat f(y,u)d\mu_{\cal F}(y)\notag
\end{align}
(this integral is viewed as a limit of finite integrals over the
complements in the leaves $F(x_i,x_j)$ of smaller and smaller
neighborhoods of the points $x_i$ $(i=1,2,3)$, and its existence
follows as above from the continuity properties of the function
$\hat f$). By our choice of $f$, for every 
$(x_1,x_2,x_3)\in BT$ the function 
$u\in G\to \nu(f)(x_1,x_2,x_3)(u)$ is contained
in $\ell^p(G)$. More precisely, the map 
$(x_1,x_2,x_3)\in BT\to \nu(f)(x_1,x_2,x_3) \in
\ell^p(\Gamma)$ is a continuous \emph{cocycle} for the action of
$G$ on $B$, i.e. it is continuous and equivariant under the action
of $G$, it satisfies $\nu\circ\sigma=({\rm sgn}(\sigma))\nu$ for
every permutation $\sigma$ of the three variables and the cocycle
identity
\begin{equation}
\nu(f)(x_2,x_3,x_4)-\nu(f)(x_1,x_3,x_4)+\nu(f)(x_1,x_2,x_4)-
\nu(f)(x_1,x_2,x_3)=0.
\end{equation}
In particular, 
for any fixed point $x\in B$ we conclude as in Section
2 that the assignment $(g_1,g_2,g_3)\to \nu(f)(g_1x,g_2x,g_3x)$
$(g_i\in G)$ defines a $G$-equivariant cocycle with values in
$\ell^p(G)$ whose cohomology class coincides with $R\Theta(f)$.

Now by a result of Adams \cite{A94} (see also \cite{Ka03} for a
more precise result), if $\sigma$ is the measure class of the
measure of maximal entropy for the geodesic flow of any convex
cocompact hyperbolic manifold whose fundamental group is a free
group with 2 generators, viewed as a $G$-invariant measure class
on the Gromov boundary $B$ of $G$, then $(B,\sigma)$ is a
\emph{strong boundary} for $G$. This means that the action of $G$
on $(B,\sigma)$ is amenable and \emph{doubly ergodic} with respect
to any separable Banach coefficient module, i.e. for every
separable Banach $G$-space $E$, every measurable $G$-equivariant
map $(B\times B,\sigma\times \sigma)\to E$ is constant almost
everywhere. As a consequence, every continuous $G$-equivariant
cocycle $BT\to \ell^p(G)$ which does not vanish identically
defines a \emph{non-vanishing} class in $H_b^2(G,\ell^p(G))$ (see
the discussion in Section 7 of \cite{M}). Thus for every $f\in
\hat{\cal H}$ such that $\nu(f)\not=0$ the class $R\Theta(f)$ does
not vanish and hence the same is true for the class $\Theta(f)$.
In other words, to show that indeed $H_b^2(\Gamma,\ell^p(\Gamma))$
is infinite dimensional we only have to find for every $m>0$ a
collection of functions $f_i\in \hat{\cal H}$ $(1\leq i\leq m)$
such that the cocycles $\nu(f_i)$ are linearly independent.

For this recall
that by Proposition 3.2 and its proof, the subspace
of $H_b^2(G,\mathbb{R})$ defined by the cohomology
classes $\Theta_G(f)\in H_b^2(G,\mathbb{R})$ 
of the quasi-morphisms $\Psi(\alpha)$ 
where $\alpha\in {\cal H}_C$ for a suitable
choice of $C\subset T$ 
is infinite dimensional (note that now we use here
the notations from Lemma 3.1 for the map $\Psi$).
On the other hand, for every $\alpha\in {\cal H}_C$
the cohomology class $\Theta_G(\alpha)
\in H_b^2(G,\mathbb{R})$
coincides with the class defined by the
continuous $\mathbb{R}$-valued cocycle
$\nu_0(\alpha):BT\to \mathbb{R}$ given by
\begin{align}
\nu_0(\alpha)(x_1,x_2,x_3)= &
\int_{F(x_1,x_2)}\tilde \alpha(y)d\mu_{\cal F}(y)\\+
\int_{F(x_2,x_3)}\tilde \alpha(y)d\mu_{\cal F}(y) &
+\int_{F(x_3,x_1)}\tilde \alpha(y)d\mu_{\cal F}(y).\notag
\end{align}
Now let $C\subset T$ and let
$\alpha_1,\dots, \alpha_m\in {\cal H}_C$ be such that
the cocycles $\nu_0(\alpha_i)$ are linearly independent;
such functions exist by Proposition 3.2 and its proof.
For every $i\leq m$ define a function
$ f_i\in \hat{\cal H}$ by $f_i(y,e)=\alpha_i(y)$ and
$f_i(y,g)\equiv 0$ for $g\not=e$. Then
\begin{align}
\nu_0(\alpha_i)(x_1,x_2,x_3)= &
\int_{F(x_1,x_2)}\sum_{u\in G}\hat f_i(y,u)d\mu_{\cal F}(y)\\+
\int_{F(x_2,x_3)}\sum_{u\in G}\hat f_i(y,u)d\mu_{\cal F}(y) &
+\int_{F(x_3,x_1)}\sum_{u\in G}\hat f_i(y,u)d\mu_{\cal F}(y)
\end{align}
and therefore since the cocycles $\nu_0(\alpha_i)$ are linearly
independent the same is true for the cocycles $\nu(f_i)$.
As a consequence, the
kernel of the map $H_b^2(\Gamma,\ell^p(\Gamma))\to
H^2(\Gamma,\ell^p(\Gamma))$ is indeed infinite dimensional.
\qed

\section{Groups acting isometrically on hyperbolic geodesic metric spaces}

In this section we consider countable groups
which admit a
weakly acylindrical isometric action on an arbitrary Gromov
hyperbolic geodesic metric space $X$. We show that the
assumptions in Theorem 3.3 are satisfied for the action of such a
group $\Gamma$ on the \emph{Gromov boundary} $\partial X$ of $X$.
From this we deduce
Theorem A from the introduction.

First recall that the Gromov boundary of a hyperbolic geodesic metric
space $X$ is defined as follows. For a fixed point $x_0\in X$,
define the \emph{Gromov product} $(y,z)_{x_0}$
based at $x_0$ of two
points $y,z\in X$  by
\begin{equation}
(y,z)_{x_0}=\frac{1}{2}\bigl(d(y,x_0)+d(z,x_0)-d(y,z)\bigr).
\end{equation}
Call two sequences $(y_i),(z_j)\subset X$ \emph{equivalent} if
$(y_i,z_i)_{x_0}\to \infty$ $(i\to \infty)$. By hyperbolicity of $X$,
this notion of equivalence defines an equivalence relation for the
collection of all sequences $(y_i)\subset X$ with the additional
property that $(y_i,y_j)_{x_0}\to \infty$ $(i,j\to\infty)$ \cite{BH}.
The boundary $\partial X$ of $X$ is the set of equivalence classes
of this relation.

The Gromov product $(\,,\,)_{x_0}$ for pairs of points in $X$
can be extended to a product on
$\partial X$ by defining \begin{equation} (\xi,\eta)_{x_0}= \sup
\liminf_{i,j\to\infty}(y_i,z_j)_{x_0}
\end{equation} where the supremum is taken over all sequences
$(y_i),(z_j)\subset X$
whose equivalence classes define the
points $\xi,\eta\in
\partial X$. For a suitable number $\chi>0$
only depending on the hyperbolicity constant of $X$ there is a
distance $\delta=\delta_{x_0}$ of bounded diameter on $\partial X$
with the property that the distance $\delta(\xi,\eta)$ between two
points $\xi,\eta\in \partial X$ is comparable to
$e^{-\chi(\xi,\eta)_{x_0}}$ (see 7.3 of \cite{GH}). More precisely,
there is a constant $\theta>0$ such that
\begin{equation}\label{dist}
e^{-\chi\theta} e^{-\chi(\xi,\eta)_{x_0}}\leq \delta(\xi,\eta)\leq
e^{-\chi(\xi,\eta)_{x_0}}\end{equation} for all $\xi,\eta\in
\partial X$. In the sequel we always assume that $\partial X$ is
equipped with such a distance $\delta$.

There is a natural topology on $X\cup \partial X$ which restricts
to the given topology on $X$ and to the topology on $\partial X$
induced by the metric $\delta$. With respect to this topology, a
sequence $(y_i)\subset X$ converges to $\xi\in
\partial X$ if and only if we have $(y_i,y_j)_{x_0}\to \infty$ and the
equivalence class of $(y_i)$ defines $\xi$. If $X$ is proper, then
$X\cup \partial X$ is compact. Every isometry of $X$ acts
naturally on $X\cup
\partial X$ as a homeomorphism. We denote by ${\rm Iso}(X)$ the
isometry group of $X$.

Since we do \emph{not} assume that $X$ is proper, for a given
pair of distinct points $\xi\not=\eta\in \partial X$ there may
not exist a geodesic $\gamma$ in $X$ connecting $\xi$ to
$\eta$, i.e. such that $\gamma(t)$ converges to $\xi$ as
$t\to-\infty$ and that $\gamma(t)$ converges
to $\eta$ as $t\to \infty$. However, there is a number
$L>1$ only depending on the hyperbolicity constant
for $X$ such that any two points in $\partial X$ can
be connected by an \emph{$L$-quasi-geodesic}.
Recall that for
$L\geq 1$, an \emph{$L$-quasi-geodesic} in $X$ is a map
$\gamma:(a,b)\to X$ for $-\infty\leq a<b\leq \infty$
such that
\begin{equation}
-L+\vert s-t\vert/L\leq d(\gamma(s),\gamma(t))\leq
L\vert s-t\vert +L
\end{equation}
for all $s,t\in (a,b)$. Note that an $L$-quasi-geodesic $\gamma$
need not be continuous. However, from every $L$-quasi-geodesic
$\gamma$ we can construct a continuous $4L$-quasi-geodesic $\tilde
\gamma$ whose \emph{Hausdorff distance}
to $\gamma$ is bounded from above
by $4L$ by replacing for each $i\geq 0$ the arc $\gamma[i,i+1]$ by
a geodesic arc $\tilde \gamma[i,i+1]$ with the same endpoints. In
other words, via changing our constant $L$ we may
assume that for any two distinct points $\xi\not=\eta \in\partial X$
there is a continuous $L$-quasi-geodesic
$\gamma$ connecting $\xi$ to $\eta$;
we then write $\gamma(-\infty)=\xi,\gamma(\infty)=\eta$
(see \cite{GH} 5.25 and 7.6; compare
also the discussion in \cite{H04}).

Recall from Section 3 the definition of a weakly hyperbolic
action of a group $G$ on a metric space of bounded diameter.
We show.

\bigskip

{\bf Lemma 4.1:} {\it Let $X$ be an arbitrary hyperbolic geodesic
metric space. Then the action of the isometry group ${\rm Iso}(X)$
on $\partial X$ is weakly hyperbolic.}

{\it Proof:} The boundary $\partial X$ of a hyperbolic geodesic
metric space is a metric space of bounded diameter where the
metric $\delta$ is constructed from the Gromov product $(\,
,\,)_{x_0}$ at a fixed point $x_0\in X$. There are numbers
$\chi>0, \theta>0$ such that inequality (\ref{dist}) above holds for our
distance $\delta$.

Our goal is to show that for every $\nu>0$ there is a constant
$\Theta=\Theta(\nu)>0$ with the following property. Let $a,b\in
\partial X$ with $\delta(a,b)\geq 2\nu$. Let $g\in {\rm Iso}(X)$
be such that $\delta(ga,gb)\geq 2\nu$; if $v\in \partial
X-\{a,b\}$ is such that $\min\{\delta(g v,g a),\delta(g v,g
b)\}\geq \nu$ then $\delta(g w,gb)\leq \Theta \delta(v,a)$ for
every $w\in \partial X$ with $\delta(w,b)\leq \nu$. Note that
since the diameter of $\partial X$ is finite, this inequality is
automatically satisfied for a suitable choice of $\Theta$ whenever
$\delta(v,a)$ is bounded from below by a universal constant. Thus
it is enough to show the claim under the additional assumption
that $\delta(v,a)\leq \epsilon$ for some fixed $\epsilon >0$
which will be determined later on.

Let $T\subset (\partial X)^3$ be the set of all triples
of pairwise distinct points in $\partial X$.
A triple $(a,b,c)\in T$
determines
(non-uniquely) an ideal $L$-quasi-geodesic triangle with vertices
$a,b,c$. The
Hausdorff distance between any two such $L$-quasi-geodesic
triangles
with the same vertices in $\partial X$ is bounded by
a universal constant.
There is a number $p_0>0$ such that for every $p\geq p_0$ and
every triple $(a,b,c)\in T$ the closed set $K(a,b,c;p)\subset X$
of all points in $X$ whose distance to each side of an
$L$-quasi-geodesic triangle with vertices $a,b,c$ is at most $p$
is non-empty. The diameter of this set is uniformly bounded by a
constant only depending on $p$ and the hyperbolicity constant for
$X$.

By the definition of the Gromov product and hyperbolicity, there
is a number $m_1>0$ with the following property. Let
$(a,b,c)\in T$ and let $\zeta$ be a continuous $L$-quasi-geodesic
connecting $b$ to $a$. Then
$\min\{(a,c)_{\zeta(0)}, (b,c)_{\zeta(0)}\}\leq m_1$ and if
$(b,c)_{\zeta(0)}\leq (a,c)_{\zeta(0)}$ then we have
$\zeta(\tau)\in K(a,b,c;m_1)$ for every $\tau\geq 0$ such that
$d(\zeta(0),\zeta(\tau))= (a,c)_{\zeta(0)}$.

Now let $\nu\in (0,1)$ and let $a,b\in \partial X$ be such that
$\delta(a,b)\geq 2\nu$. By
hyperbolicity and inequality (\ref{dist}) above, there is a
constant $m_0=m_0(\nu)>0$ such that every continuous
$L$-quasi-geodesic connecting two points $a\not=b\in \partial X$
with $\delta(a,b)\geq \nu$ intersects the ball $B(x_0,m_0)$. Let
$\gamma$ be a continuous $L$-quasi-geodesic connecting
$b=\gamma(-\infty)$ to $a=\gamma(\infty)$ which is parametrized in
such a way that $\gamma(0)\in B(x_0,m_0)$. Let $\theta >0,\chi>0$
be as in inequality (\ref{dist}), let $R_0=
\chi(m_0+m_1+\theta)$ and let $v\in
\partial X-\{a,b\}$ be such that $\delta(a,v)\leq e^{-R_0}$; then
$\delta(a,v)=e^{-R}$ for some $R\geq R_0$. By inequality
(\ref{dist}) we have $R/\chi-\theta\leq (a,v)_{x_0}\leq R/\chi$
and hence 
\begin{equation}
R/\chi-\theta-m_0\leq (a,v)_{\gamma(0)}\leq R/\chi+m_0
\end{equation}
since $d(x_0,\gamma(0))\leq m_0$. From the assumption on $R$ we
obtain that $(a,v)_{\gamma(0)}\geq m_1$ and hence
$\gamma(\tau)\in K(a,b,v;m_1)$ for all $\tau
\geq 0$ such that $d(\gamma(0),\gamma(\tau))=
(a,v)_{\gamma(0)}$.

Let $g\in {\rm Iso}(X)$ be such that $\delta(ga,gb)\geq 2\nu$ and
$\min\{\delta(ga,gv),\delta(gb,gv)\}\geq \nu$. Then the
$L$-quasi-geodesic $g\gamma$ intersects $B(x_0,m_0)$ and the same
if true for any $L$-quasi-geodesic connecting $ga$ to $gv$ or
connecting $gb$ to $gv$ and consequently $x_0\in K(ga,gb,gv;m_0)$.
If as before $\tau>0$ is such that
$d(\gamma(0),\gamma(\tau))=(a,v)_{\gamma(0)}$ 
then $\gamma(\tau)\in
K(a,b,v;m_1)$ and therefore 
\begin{equation}
\{x_0,g\gamma(\tau)\}\subset
K(ga,gb,gv;m_0+m_1) =gK(a,b,v;m_0+m_1). 
\end{equation}
Now the
diameter of the set $K(ga,gb,gv;m_0+m_1)$ is bounded from
above by a constant $m_2=m_2(\nu)>0$ only depending on $\nu$ and
hence $d(g\gamma(\tau),x_0)\leq m_2$.

Let $w\in \partial X$ be such that $\delta(w,b)\leq \nu$. Then
$\delta(w,a)\geq \nu$ and by inequality (\ref{dist}) above, the
Gromov product $(w,a)_{x_0}$ is bounded from above by a universal
constant and the same is true for $(w,a)_{\gamma(0)}$. In
particular, the $L$-quasi-geodesic ray $\gamma[0,\infty)$
connecting $\gamma(0)$ to $a$ is contained
in a uniformly bounded neighborhood of any $L$-quasi-geodesic
connecting $w$ to $a$. With $\tau>0$ as above
we have $\vert d(\gamma(\tau),\gamma(0))-R/\chi\vert 
\leq m_0+\theta$ and hence
by the definition of the Gromov product and hyperbolicity,
the quantity $(b,w)_{\gamma(\tau)}-R/\chi=
(gb,gw)_{g\gamma(\tau)}- R/\chi$ is bounded from below by a
universal constant. But $d(g\gamma(\tau),x_0)\leq m_2$ and hence
we have 
\begin{equation}
\vert (gb,gw)_{g\gamma(\tau)}-(gb,gw)_{x_0}\vert =\vert
(b,w)_{\gamma(\tau)}-(gb,gw)_{x_0}\vert \leq m_2. 
\end{equation}
Using once more
the estimate (\ref{dist}) we conclude that there is a number
$\Theta>1$ only depending on $\nu$ such that $\delta(gb,gw)\leq
\Theta e^{-R}=\Theta\delta(a,v)$. This shows that the action of
${\rm Iso}(X)$ on $\partial X$ is weakly hyperbolic. \qed

\bigskip

As in the introduction, we call an isometric action on $X$ of a
countable group $\Gamma$ \emph{weakly acylindrical} if for every
point $x_0\in X$ and every $m>0$ there are numbers $R(x_0,m)>0$
and $c(x_0,m)>0$ with the following property. If $x,y\in X$ with
$d(x,y)\geq R(x_0,m)$ are such that a geodesic $\gamma$ connecting
$x$ to $y$ meets the $m$-neighborhood of $x_0$ then there are at
most $c(x_0,m)$ elements $g\in \Gamma$ such that $d(x,gx)\leq m$
and $d(y,gy)\leq m$. We have.

\bigskip

{\bf Lemma 4.2:} {\it Let $X$ be a hyperbolic geodesic metric
space and let $\Gamma$ be a countable subgroup of ${\rm Iso}(X)$
whose action on $X$ is weakly acylindrical; then the action of
$\Gamma$ on the space or triples of pairwise distinct points in
$\partial X$ is metrically proper.}

{\it Proof:} Let $X$ be a hyperbolic geodesic metric space and let
$\Gamma$ be a countable subgroup of ${\rm Iso}(X)$ whose action on
$X$ is weakly acylindrical.
Then $\Gamma$ acts as a
group of homeomorphisms on the Gromov boundary $\partial X$ of
$X$. Recall that $\partial X$ is a metric space of bounded
diameter where the metric $\delta$ is constructed from the Gromov
product $(\, ,\,)_{x_0}$ at a fixed point $x_0\in X$ and it
satisfies the estimate (\ref{dist}) from the beginning of this
section for some $\chi >0,\theta >0$ and all $\xi\not=\eta\in
\partial X$. We have to show that the action of $\Gamma$ on the
space $T$ of triples of pairwise distinct points in $\partial X$
is metrically proper.

For this let $\nu>0$ be fixed. There are numbers $L\geq 1,
m_0=m_0(\nu)>0$ such that any two points $x\not=y \in \partial X$
can be connected by a continuous $L$-quasi-geodesic, and if
$\delta(x,y)\geq \nu$ then this quasi-geodesic intersects the ball
$B(x_0,m_0)$.

By hyperbolicity, the Hausdorff distance between
any two $L$-quasi-geodesics connecting the same
points in $\partial X$ is bounded from above by a universal
constant. Moreover, there is a universal constant 
$m_1=m_1(\nu)>m_0$
with the following property. Let $a\not=b, x\not=y
\in \partial X$ and assume that $\delta(a,b)\geq 2\nu$ and that
for some $R>-\log \nu/2$ we have $\delta(a,x)\leq e^{-R},
\delta(b,y)\leq e^{-R}$. Let $\gamma$
be a continuous $L$-quasi-geodesic connecting $b=\gamma(-\infty)$
to $a=\gamma(\infty)$ 
and let $\eta$ be a continuous $L$-quasi-geodesic connecting
$y=\eta(-\infty)$ to $x=\eta(\infty)$;
then $\gamma, \eta$ intersect the ball $B(x_0,m_0)$,
and the
intersection of $\gamma$
with the ball $B(x_0,R/\chi)$ is contained in
the $m_1$-neighborhood of $\eta$. 

As in the proof of Lemma 4.1, 
for $p>0$ and a triple $(u,v,w)\in
T$ of pairwise distinct points in $\partial X$ let
$K(u,v,w;p)\subset X$ be the set of all points whose distance to
each side of an $L$-quasi-geodesic triangle with vertices $u,v,w$
is at most $p$. 
By the arguments in the proof of Lemma 4.1
there is a constant $m_2>m_1$ with the following
property. Let $x,y\in \partial X$ with
$d(x,y)\geq\nu$. If $z\in \partial X$ and $k\geq -\log \nu/2$ 
are such that $e^{-k}\leq \delta(x,z)\leq e^{-k+1}$
then the distance between
$x_0$ and $K(x,y,z;m_0)$ is contained in the
interval $[k/\chi-m_2,k/\chi+m_2]$.
The diameter of the sets
$K(x^\prime,y^\prime,z^\prime;m_0)$ is bounded
from above by a universal constant $\rho>0$ 
only depending on $m_0$ and the hyperbolicity
constant of $X$.

Let $(a,b,c)\in (\partial X)^3$ be a 
triple of points whose pairwise distance
is at least $2\nu$. Let $R\geq -\log \nu/2$
be a number to be determined later, let $U(a),U(b),U(c)$
be the open $e^{-R}$-neighborhood of $a,b,c$ in $\partial X$ and
let $x\in U(a),y\in U(b),z\in U(c)$. Let $a^\prime,b^\prime\in
\partial X$ be such that $\delta(a^\prime,b^\prime)\geq 2\nu$ and
assume that there is some $g\in \Gamma$ such that $gx=a^\prime,
gy=b^\prime$ and 
$\delta(gz,a^\prime)\in [e^{-k},e^{-k+1}]$ for some $k\geq R$.
Then $g$ maps a continuous $L$-quasi-geodesic $\eta$ 
connecting $y$ to $x$ 
with $\eta(0)\in K(x,y,z;m_0)$
to a continuous
$L$-quasi-geodesic $g\eta$ connecting $b^\prime$ to 
$a^\prime$. Since $g(\eta(0))\in K(a^\prime,b^\prime,gz;m_0)$
we have 
\begin{equation}\label{distance}
\vert d(g\eta(0),x_0)-k/\chi\vert \leq m_2+\rho.
\end{equation}

Now let $x^\prime\in U(a),y^\prime\in U(b),
z^\prime\in U(c)$ and let $g^\prime\in \Gamma$
be such that $g^\prime x^\prime=a^\prime=gx,
g^\prime y^\prime=b^\prime=gy$ and 
$\delta(g^\prime z^\prime,a^\prime)\in [e^{-k},e^{-k+1}]$.
Let $\eta^\prime$ be continuous $L$-quasi-geodesics
connecting $y^\prime$ to $x^\prime$, 
with $\eta^\prime(0)\in B(x_0,m_0)$.
As above, let $\gamma$
be a continous $L$-quasi-geodesic
connecting $b$ to $a$ with
$\gamma(0)\in B(x_0,m_0)$ and let $\sigma<0$ be such that
$d(x_0,\gamma(\sigma))=R/\chi$.
Then there are numbers $\tau<0,\tau^\prime<0$ such that
$d(\eta(\tau),\gamma(\sigma))\leq m_1,
d(\eta^\prime(\tau^\prime),\gamma(\sigma))\leq m_1$
and therefore
$d(\eta(\tau),\eta^\prime(\tau^\prime))\leq 2m_1.$
In particular, we have 
\begin{equation}
\vert d(\eta(0),\eta(\tau))-d(\eta^\prime(0),
\eta^\prime(\tau^\prime))\vert \leq 2m_0+2m_1
\end{equation}
The images of $\eta,\eta^\prime$ under $g,g^\prime$ 
are continuous
$L$-quasi-geodesics connecting
$b^\prime$ to $a^\prime$. The estimate (\ref{distance})
is valid for $g^\prime$ as well
and hence by hyperbolicity, the distances $d(g\eta(0),
g^\prime(\eta^\prime(0))),$  
$d(g(\eta(\tau)),g^\prime(\eta^\prime(\tau^\prime)))$ are  
bounded from above by a universal constant
$m_3>2m_2$. 
Together we conclude that
\begin{equation}
d(g^{-1}g^\prime(\eta(0)),\eta(0))\leq 2m_3,\quad
d(g^{-1}g^\prime(\eta(\tau)),
\eta(\tau))\leq 2m_3.
\end{equation}

Now if $R_0=R(x_0,2m_3)$ is as in the definition of a weakly
acylindrical action, then for $R\geq \chi R_0$ and any $k\geq R$
the number of elements $g,g^\prime\in \Gamma$ with this property is
bounded from above by for a universal constant independent of $R$
and $k$. This shows that the action of $\Gamma$ on $\partial X$
satisfies the first property in the definition of a metrically
proper action.

The second property in the definition of a metrically proper
action follows from exactly the same argument. Namely, using our
above notation, there is a number $\kappa >m_0(\nu)$ only
depending on $\nu$ such that if $Z\subset\partial X$ is the set of
all points whose distance to $U(a),U(b)$ is at least $\nu$ then
there is a number $\tau_0>0$ such that for any $x\in U(a),y\in
U(b)$ and $z\in Z$ the set $K(x,y,z)$ is contained in the ball of
radius $\kappa>0$ about $x_0$. In other words, for any element
$g,h\in \Gamma$ which map a triple $(x,y,z)\in U(a)\times
U(b)\times Z$ to a triple of points whose pairwise distance is
bounded from below by $\nu$, the distance between $x_0$ and $gx_0$
is at most $\kappa$. Our above consideration then shows that we
can find a number $\tilde R(\nu)>0$ 
depending on $\nu$ and some $\tilde m(\nu)>0$
such that the second
requirement in the definition of a metrically proper action
holds with these constants
and for the action of $\Gamma$ on $\partial X$.\qed

\bigskip

Recall from Section 3 the definition of a homeomorphism
with north-south dynamics
of a metric space of finite diameter.
The statement of the next simple
lemma is well known in the case that the hyperbolic
space $X$ is proper;
we include a short proof for the sake of completeness
since we did not find
a suitable reference for the general case.

\bigskip

{\bf Lemma 4.3:} {\it Let $X$ be a hyperbolic geodesic
metric space and let $g$ be an isometry
of $X$ such that for some $x\in X$ the map
$k\to g^kx$ is a quasi-isometric embedding of the integers
into $X$. Then $g$ acts on $\partial X$ with north-south
dynamics.}

{\it Proof:} Let $g$ be an isometry of the hyperbolic
geodesic metric space $X$ with the property that
for some $x\in X$ the map $k\to g^k x$ is a quasi-isometric
embedding of the integers into $X$. Then
the sequence
$(g^kx)_{k\geq 0}\subset X$ converges to a point
$a\in \partial X$,
and the sequence $(g^{-k}x)_{k\geq 0}\subset X$ converges
to a point $b\in \partial X-\{a\}$.
The limit set of the infinite cyclic
group $G$ generated by $g$
consists of the two points $a\not= b\in \partial X$,
and these are fixed points
for the action of $G$ on $\partial X$.

By hyperbolicity
there is a number $m>0$ such that
for every $\xi\in \partial X-\{a,b\}$ the closed
set $K(a,b,\xi;m)\subset X$ of all points in $X$
whose distance
to each side of an $L$-quasi-geodesic
triangle with vertices $a,b,\xi$ is at most $m$
is non-empty and its
diameter $K(a,b,\xi;m)$
is bounded independently of $\xi$.
Since the assignment
$k\to g^k(x)$ is a quasi-isometric embedding of the integers
into $X$, we may assume by possibly enlarging $m$
that each of the sets $K(a,b,\xi;m)$
intersects $Q=\{g^k(x)\mid
k\in \mathbb{Z}\}$. Thus
there is a number $\ell >0$ and for every
$\xi\in \partial X-\{a,b\}$ there is some $\kappa(\xi)\in\mathbb{Z}$ such that
the set $\{g^\kappa(x)\mid \kappa(\xi)\leq \kappa\leq
\kappa(\xi)+\ell\}$ contains the intersection of $K(a,b,\xi;m)$ with
$Q$.
Then $\vert \kappa(g^j\xi)-\kappa(\xi)-j\vert \leq \ell$
for all $j\in \mathbb{Z}$
and hence the set $D=\{\xi\in \partial X-\{a,b\}\mid
0\leq \kappa(\xi)\leq \ell\}$  does not contain $a,b$ in its
closure and it satisfies
$\cup_{j\in \mathbb{Z}}g^jD=\partial X-\{a,b\}$.
Moroever, for every neighborhood $U$ of $a$,
$V$ of $b$ there is a number
$j>0$ such that $g^j(X-V)\subset U,g^{-j}(X-U)\subset V$.
Hence the isometry
$g$ acts with north-south dynamics on $\partial X$.
This shows the lemma.
\qed

\bigskip

Call an isometry of $X$ \emph{hyperbolic} if it acts on
$\partial X$ with north-south dynamics with respect to some fixed
points $a\not= b$. The following 
corollary is immediate
from Lemma 4.1, Lemma 4.2 and the remark after
Proposition 3.2 in Section 3. We refer to \cite{PR04}
for a similar result for the group $SL(2,\mathbb{Z})$.

\bigskip

{\bf Corollary 4.4:} {\it Let $\Gamma$ 
be a countable group which
admits a weakly acylindrical isometric action
on a hyperbolic geodesic metric
space. Let $g_1,\dots,g_k\in \Gamma$ be 
hyperbolic elements with ordered pairs of fixed
points $(a_i,b_i)$. If the $\Gamma$-orbits of 
$(a_i,b_i),(b_i,a_i)$ are pairwise disjoint then
for every $(q_1,\dots,q_k)\in \mathbb{R}^k$
there is a quasi-morphism $\phi$ for $\Gamma$
with $\lim_{\ell\to\infty}\phi(g_i^\ell)/\ell=q_i$
for every $i\leq k$.}

\bigskip

The \emph{limit set} of an isometric action of a group $\Gamma$ on
$X$ is the set of accumulation points in $\partial X$ of an orbit
$\Gamma x$ $(x\in X)$ of $\Gamma$; it does not depend on the
orbit. A subgroup $\Gamma$ of ${\rm Iso}(X)$ is called
\emph{elementary} if its limit set contains at most 2 points.
The next result is Theorem A from the introduction.

\bigskip

{\bf Theorem 4.5:} {\it Let $\Gamma$ be a countable group which
admits a non-elementary weakly acylindrical isometric action on a
Gromov hyperbolic geodesic metric space $X$; then the kernels of
the natural homomorphisms $H_b^2(\Gamma,\mathbb{R})\to
H^2(\Gamma,\mathbb{R}), H_b^2(\Gamma, \ell^p(\Gamma))\to
H^2(\Gamma, \ell^p(\Gamma))$ $(1< p<\infty)$ are infinite
dimensional.}

{\it Proof:} Let $X$ be a hyperbolic geodesic metric space and let
$\Gamma$ be a countable non-elementary weakly acylindrical
subgroup of the isometry group of $X$. By assumption, the limit
set $\Lambda$ of $\Gamma$ contains at least 3 points.
Then this limit set is a
$\Gamma$-invariant closed subset of $\partial X$
without isolated points (see \cite{GH}). Our goal is to
show that the action of $\Gamma$ on $\Lambda$ satisfies the
assumptions in Theorem 3.3.

By Lemma 4.1 and Lemma 4.2, the action of 
$\Gamma$ on $\Lambda$ is
weakly hyperbolic and the action of $\Gamma$ on the space of
triples of pairwise distinct points in $\Lambda$ is metrically
proper. Using Lemma 4.3 it is enough 
to show that $\Gamma$ contains a free
subgroup $G$ with two generators
which has the following additional properties.
\begin{enumerate}
\item For some $x\in X$ the orbit map $g\in G\to gx\in X$
is a quasi-isometric embedding of $G$ into $X$.
\item
There are infinitely many $g_i\in G$ $(i\geq 0)$ such that the
ordered pairs of fixed points of $g_i,g_j^{-1}$ 
are contained in pairwise distinct orbits of
the action of $\Gamma$ on $\Lambda\times \Lambda$.
\end{enumerate}
Note that the first property guarantees that there
is a continuous $G$-equivariant embedding of the
Gromov boundary $B$ of $G$ into $\Lambda$.

The existence of a free group $G$ with two generators
and with property (1) above is immediate from
the ping-pong lemma and our requirement that
the group $\Gamma$ is non-elementary (compare 
\cite{GH}).

Now let $e\not= g\in G$ and let 
$(a,b)$ be
the ordered pair of fixed points of the action of $g$ on
$\partial X$. Choose a closed 
subset of $\partial X$ which is contained in
$X-\{a,b\}$ and is a fundamental domain
$D$ for the action on $\partial X-\{a,b\}$ of the
infinite cyclic subgroup of $G$ generated by $g$.
Assume that there is a sequence $(a_i,b_i)\in 
\partial X\times \partial X$
contained in the $\Gamma$-orbit of $(a,b)$ 
with $(a_i,b_i)\to (a,b)$.
Let $\delta$ be a Gromov distance on $\partial X$
and write $\nu=\min\{\delta(a,b),\delta(\{a,b\},D)\}/4$. 
Let $R(\nu)>0$ be as in the definition of a metrically
proper action for $\Gamma$ and let
$U,V$ be the open $e^{-R(\nu)}$-neighborhood
of $a,b$. For sufficiently large $i$ we have $a_i\in U,
b_i\in V$. By our assumption, 
there are $h_i\in \Gamma$ such that
$h_ia_i=a,h_ib_i=b$. Then $h_i^{-1}gh_i$ is a hyperbolic
isometry with fixed points $a_i,b_i$.
Since a hyperbolic
isometry fixes \emph{precisely} two points in $\partial X$,
the elements $h_i$ are pairwise distinct and the same is
true for their compositions with an arbitrary
power of $g$. Namely, otherwise there are $i\not=j$
and $\ell\in \mathbb{Z}$ such that
$g^\ell=h_ih_j^{-1}$ which contradicts the fact
that $(a,b)$ are fixed points for $g$, 
$(a_i,b_i)\not=(a_j,b_j)$ and that $h_i$ is a homeomorphism.
However, by the choice of $D$ there
is for each $i>0$ some $k(i)\in \mathbb{Z}$ such that
$g^{k(i)}h_i D\cap D\not=\emptyset$ and hence
$g^{k(i)}h_i(U\times V\times D)\cap
U\times V\times D\not=\emptyset$ for all sufficiently
large $i$. This contradicts
our assumption that the action of $\Gamma$ on 
the space of triples of pairwise
distinct points in $\partial X$ is metrically proper.

As a consequence,
for every ordered pair $(a,b)$ of fixed points
of an element $e\not=g\in G$ the $\Gamma$-orbit
of $(a,b)$ is a \emph{discrete} subset of 
$\partial X\times \partial X-\Delta$ (note that
this fact has already been established in
the proof of Proposition 3.2). 
Since on the other hand the sets of pairs of fixed
points for the elements of $G$ are \emph{dense}
in $B\times B-\Delta$, there are infinitely many
such pairs $(a_i,b_i)$
which are pairwise contained in
distinct orbits under the action of $\Gamma$.
Our argument also implies that we may in addition require
that the ordered pairs $(a_i,b_i)$ are not
contained in the $\Gamma$-orbit of $(b_j,a_j)$
for any $j$. 
 
We use this fact to show that
we can find infinitely many $g_i\in G$ with
the property that the $\Gamma$-orbits of the
ordered pairs of fixed points $(a_i,b_i),
(b_j,a_j)$ of $g_i,g_j^{-1}$ are all disjoint
(see the argument in \cite{BF}).
Namely, choose two independent elements
$g_1,g_2\in G$ which generate a free subgroup
with the property that the ordered
pairs of fixed points $(a_1,b_1),(b_1,a_1)$ of 
$g_1,g_1^{-1}$
are not contained in the $\Gamma$-orbit of 
the ordered pair of fixed points
$(b_2,a_2)$ of $g_2^{-1}$. We may assume that
the group generated by $g_1,g_2$ equals $G$ and that
there is an $L$-quasi-isometric $G$-equivariant embedding
$\rho$ of the Cayley graph $CG$ of $G$ into $X$ which
induces an equivariant embedding of the
Gromov boundary $B$ of $G$ into $\partial X$.
Identify $CG$ with its image under our embedding. 
For $0<<n_1<<m_1<<n_2<<m_2$ consider the
element $f=g_1^{n_1}g_2^{m_1}g_1^{n_2}g_2^{m_2}\in G$.
If $\gamma$ is the axis of $f$ in $CG$ and if
$h\in \Gamma$ maps the ordered pair $(a,b)$ of
fixed points for $f$ to $(b,a)$, then it maps the
inverse $\rho(\gamma)^{-1}$ of $\rho(\gamma)$ into
a uniformly bounded neighborhood of $\rho(\gamma)$.
Now a fundamental domain for the action of $f$ on its
axis $\gamma$ is composed of four arcs $\gamma_1,\dots,\gamma_4$
where $\gamma_1$ is the geodesic arc in $CG$ connecting
$e$ to $g_1^{n_1}$, $\gamma_2$ is the translate
under $g_1^{n_1}$ of the geodesic arc connecting $e$
to $g_2^{m_1}$ etc. 
As a consequence, there is a subsegment of
the axis of a conjugate of $g_1$ in $G$ 
whose length tends to infinity as $n_1\to \infty$ and 
which is mapped by $h$ into a uniformly bounded
neigborhood of a 
subsegment of the axis of a conjugate of $g_2^{-1}$
(see \cite{BF}).
For sufficiently large $n_1$ 
this violates our observation that the
$\Gamma$-orbits of $(a_i,b_1),(b_2,a_2)$ are
discrete and disjoint.

As a consequence, property (2) above holds for
$G$ as well (compare also the discussion in 
\cite{BF}). 
Thus our theorem is a consequence of Theorem 3.1.
\qed

\section{Applications}

This section is devoted to a discussion of applications
of our Theorem A from the introduction.
We begin with the proof of 
Corollary B from the
introduction. For this let $S$ be an oriented surface of genus
$g\geq 0$ with $m\geq 0$ punctures. We assume that $S$ is
\emph{non-exceptional}, i.e. that $3g-3+m\geq 2$. The
\emph{complex of curves} ${\cal C}(S)$ for $S$ is the simplicial
complex whose vertices are free homotopy classes of
\emph{essential simple closed curves} on $S$, i.e.
simple closed curves which are neither contractible nor
freely homotopic into a puncture of $S$. The
simplices in ${\cal C}(S)$ are spanned by collections of such
curves which can be realized disjointly. Since $S$ is
non-exceptional by assumption, the complex of curves is connected.
If we equip each simplex in ${\cal C}(S)$ with the standard
euclidean metric of side-length one, then we obtain a length
metric on  ${\cal C}(S)$, and this length metric defines on ${\cal
C}(S)$ the structure of a hyperbolic geodesic metric space.
However, ${\cal C}(S)$ is not locally finite and hence this
geodesic metric space is not locally compact
(for all this see \cite{MM,B02,H05}). A
description of its Gromov boundary is contained in \cite{K99,H04}.

The \emph{mapping class group} ${\cal M}_{g,m}$ of $S$ is the
group of isotopy classes of orientation preserving homeomorphisms
of $S$. It acts as a group of isometries on the complex of curves
${\cal C}(S)$ of $S$. Bowditch \cite{B03} showed that this action
is weakly acylindrical. Thus we can apply Theorem 4.5 and
deduce Corollary B from the
introduction which extends the result of Bestvina
and Fujiwara \cite{BF}.

\bigskip

{\bf Proposition 5.1:} {\it Let $\Gamma$ be an arbitrary subgroup of
${\cal M}_{g,m}$ which is not virtually abelian;
then the group $H_b^2(\Gamma,\mathbb{R})$
is infinite dimensional. If moreover $\Gamma$ does not
contain a normal subgroup which 
virtually splits as a direct product of two infinite
groups then for every $(1< p<\infty)$
the group $H_b^2(\Gamma,\ell^p(\Gamma))$  is
infinite dimensional as well.}

{\it Proof:} Recall from \cite{MP} the classification
of subgroups $\Gamma$ of ${\cal M}_{g,m}$. There are 4 cases.
\begin{enumerate}
\item $\Gamma$ contains two independent \emph{pseudo-Anosov} elements.
\item The limit set of the action of
$\Gamma$ on ${\cal C}(S)$ consists of precisely two points $a\not=b$.
\item $\Gamma$ is finite.
\item $\Gamma$ preserves a nontrivial
system of pairwise disjoint essential simple closed mutually not freely
homotopic curves on $S$.
\end{enumerate}

The action of the mapping
class group on ${\cal C}(S)$ is weakly acylindrical
\cite{B03} and
hence the same is true for the action of an
arbitrary subgroup $\Gamma$ of ${\cal M}_{g,m}$.
If $\Gamma$ is as in case (1) above then the limit set of
$\Gamma$ contains at least 3 points and therefore
$\Gamma$
is a non-elementary subgroup
of the isometry group of ${\cal C}(S)$. By Theorem 4.5,
the groups $H_b^2(\Gamma,\mathbb{R}), H_b^2(\Gamma,\ell^p(\Gamma))$
are infinite dimensional for every $p\in (1,\infty)$.

In the case (2) above, each element of $\Gamma$ maps a quasi-geodesic
connecting $a$ to $b$ into a uniformly bounded neighborhood of
itself. Since the action of $\Gamma$ on ${\cal C}(S)$ is weakly
acylindrical, the group $\Gamma$ is virtually cyclic
(compare the discussion in \cite{BF}).

In the case (4) there is a maximal system ${\cal S}$ of pairwise
disjoint essential simple closed mutually not freely homotopic
curves preserved by $\Gamma$.
If we cut $S$ open along ${\cal S}$
and replace each boundary
circle of the resulting bordered surface by a puncture then
we obtain a possibly
disconnected surface $S^\prime$ of finite type and
of bigger Euler characteristic. There is a natural
homomorphism of $\Gamma$ onto
a subgroup $\Gamma^\prime$ of the mapping class
group of $S^\prime$. Its kernel $K$ is a free abelian
group generated by Dehn twists about the curves of our curve
system. Thus by Theorem 12.4.2 of \cite{M} (see also
Corollary 3.6 of \cite{MS05}), the natural
map $H_b^2(\Gamma^\prime,\mathbb{R})\to 
H_b^2(\Gamma,\mathbb{R})$ is an isomorphism.

Let $S_1^\prime,\dots, S_p^\prime$ be the
connected components of $S^\prime$.
An element $g\in \Gamma^\prime$ permutes the components
of $S^\prime$. This means that there is a homomorphism
$\kappa$ of $\Gamma^\prime$ into the group
of permutations of $\{1,\dots,p\}$ 
whose kernel is the normal subgroup 
$G$ of $\Gamma$ of all elements which fix
each component $S_i^\prime$.
Thus there is an exact sequence
\begin{equation}
0\to G\to \Gamma^\prime\to Q\to 0
\end{equation}
where $Q$ is a finite group. This sequence induces an exact
sequence \cite{M}
\begin{equation}
\dots\to H_b^2(Q,\mathbb{R})\to H_b^2(\Gamma^\prime,\mathbb{R})
\to H_b^2(G,\mathbb{R})\to H_b^3(Q,\mathbb{R})\to\dots
\end{equation}
Since the group $Q$ is finite, its bounded cohomology
with real coefficients 
is finite dimensional and therefore we conclude that
$H_b^2(\Gamma^\prime,\mathbb{R})$ is infinite dimensional
if and only if this is the case for
$H_b^2(G,\mathbb{R})$.

For $i\leq p$ denote by $G_i$ the
projection of $G$ to a subgroup of the mapping
class group of $S_i^\prime$.
If $G_i$
preserves a non-trivial system ${\cal S}_i$ 
of pairwise disjoint essential simple closed
not mutually freely homotopic curves on $S_i^\prime$
then the $\Gamma^\prime$-translates
of this system is a $\Gamma^\prime$-invariant
curve system on $S^\prime$ which lifts
to a $\Gamma$-invariant curve system on $S$
strictly containing ${\cal S}$. This
contradicts the maximality of the system ${\cal S}$.

An exceptional component $S_i^\prime$ of $S^\prime$ either is a
thrice punctured sphere with finite mapping class group, or
$S_i^\prime$ is a once punctured torus or a forth punctured sphere
with word hyperbolic mapping class group. 
Therefore either $\Gamma^\prime$ and hence $\Gamma$
is virtually abelian
or after reordering, the group $G_1$ admits a weakly
acylindrical action as 
a non-elementary group of isometries
on a hyperbolic geodesic metric space.
In particular, if $\Gamma$ is not virtually
abelian then the second bounded cohomology
group $H_b^2(G_1,\mathbb{R})$ is infinite dimensional.

Let $R$ be the kernel of the homomorphism
$G\to G_1$. Then we have an exact sequence
\begin{equation}\label{exactsequence}
0\to R\to G\to G_1\to 0. 
\end{equation} 
Since necessarily $H_b^1(R,\mathbb{R})=0$ (see \cite{M})
we obtain from the induced exact sequence 
of bounded cohomology groups that
$H_b^2(G,\mathbb{R})$ is infinite dimensional
if this is the case for $H_b^2(G_1,\mathbb{R})$.
In other words, either $\Gamma$ is virtually abelian
or the second bounded cohomology group
$H_b^2(\Gamma,\mathbb{R})$ is infinite dimensional.

We are left with investigating the groups
$H_b^2(\Gamma,\ell^p(\Gamma))$. 
Using our above notations, note first that
if the kernel $K$ of the natural projection
$\pi:\Gamma\to \Gamma^\prime$ is nontrivial,
then the normal subgroup $\pi^{-1}(G)$ of
$\Gamma$ splits as a direct product. Thus as before,
we may assume that $\Gamma=\Gamma^\prime$.
Then $H_b^2(\Gamma,\ell^p(\Gamma))$ is infinite
dimensional if this is the case for
$H_b^2(G,\ell^p(G))$. 
Namely, if the centralizer
$Z_\Gamma(G)$ of $G$ in $\Gamma$ is infinite then 
the center of $G$ is infinite and hence either
$G$ is virtually abelian or  
$G$ splits as a direct product of two infinite
groups. Thus we may assume that $Z_\Gamma(G)$ 
is finite. Then every function $f\in \ell^p(G)$ which
is invariant under the action of the
finite center of $G$ defines a function in
the $G$-module $\ell^p(\Gamma)^{Z_\Gamma(G)}$ of
$Z_\Gamma(G)$-invariant points in $\ell^p(\Gamma)$
which vanishes outside of $GZ_\Gamma(G)$.
It follows that the second bounded
cohomology group $H_b^2(G,\ell^p(\Gamma)^{Z_\Gamma(G)})$
is infinite dimensional. The finite group
$Q$ admits an isometric action on 
$H_b^2(G,\ell^p(\Gamma)^{Z_\Gamma(G)})$
induced from the action of $Q$ on $G$ by conjugation
(Corollary 8.7.3 of \cite{M}). 
Unsing the explicit form of this action we conclude that
the subspace of $H_b^2(G,\ell^p(\Gamma)^{Z_\Gamma(G)})$
of elements
which are fixed by $Q$ is infinite dimensional if
this is the case for $H_b^2(G,\ell^p(G))$.
On the other hand, since the
group $G$ is infinite by assumption, there is
no nonzero $G$-invariant vector in
$\ell^p(\Gamma)$ and hence by the Hochschild-Serre
spectral sequence for bounded cohomology
(Theorem 12.0.3 of \cite{M}),  
the second bounded cohomology $H_b^2(\Gamma,\ell^p(\Gamma))$
is infinite dimensional if this is the case
for $H_b^2(G,\ell^p(G))$.

Let $N_j$ be the kernel of the projection
of $G$ onto a subgroup of the mapping class group of
$S^\prime-S_j^\prime$. Then $N_j$ consists of mapping
classes 
which act trivially on $S_i$ for all $i\not=j$. For $i\not=j$,
the groups $N_i,N_j$ only intersect in the identity
and commute. Thus if the groups
$N_i,N_j$ are infinite for some $i\not= j$, then $G$
contains a normal subgroup which is the direct
product of two infinite groups. The smallest
normal subgroup of $\Gamma$ containing $N_i,N_j$ 
contains the direct product of $N_i,N_j$ as
a subgroup of finite index, i.e. this normal subgroup
virtually splits as a direct product.
Thus for the purpose
of our proposition we may assume after reordering that
$N_i$ is finite for all $i>1$.

Consider first the case
that $N_1$ is infinite. Denote as before by $R$ the kernel
of the natural projection $G\to G_1$ into the mapping class group
of $S_1^\prime$. 
The subgroup of $G$ generated by $N_1,R$ is normal
and the direct product of $N_1$ and $R$. Hence as above,
if $\Gamma$ does not contain a normal subgroup 
which virtually splits as a direct product of two infinite
groups then $R$ is finite, and the
quotient group $G/R$ can naturally
be identified with the group $G_1$.

Assume that this holds true.
By Theorem 4.5 and our assumption that
$G$ is not virtually abelian, the second bounded
cohomology group $H_b^2(G_1,\ell^p(G_1))$
is infinite dimensional for every $p\in (1,\infty)$.
Now the group $R$ is finite and therefore
averaging over the orbits of the action of 
$R$ shows that $\ell^p(G_1)$ as a $G_1$-module
can naturally be identified with
the $G_1$-module $\ell^p(G)^R$ 
of all $R$-invariant points
in $\ell^p(G)$. As a consequence,
the group $H_b^2(G_1,\ell^p(G)^R)$ is infinite
dimensional, and therefore from the 
Hochschild-Serre spectral sequence
(Theorem 12.0.3 of \cite{M})
we deduce that the same is true for
$H_b^2(G,\ell^p(G))$.
As a consequence, if $N_1$ is infinite
and if $\Gamma$ does not contain a normal
subgroup which virtually splits as a direct product then
$H_b^2(G,\ell^p(G))$ is infinite
dimensional as claimed.

Finally we have to consider the case that
$N_1$ is finite, i.e. that the kernel of the natural projection
of $G$ to a subgroup of the mapping class
group of $S_2^\prime\cup\dots\cup S_p^\prime$
is finite. By our above consideration,
for every $p\in (1,\infty)$ the group
$H_b^2(G,\ell^p(G))$ is infinite
dimensional if this is the case for
$H_b^2(G/N_1,\ell^p(G/N_1))$.
Since $\Gamma$ contains a 
normal subgroup which virtually splits
as a direct product if this is the case for
$G/N_1$, an application of the above consideration
to the group $G/N_1$ yields inductively the following.
Either $\Gamma$
contains a normal subgroup which
virtually splits as a direct product or 
$H_b^2(\Gamma,\ell^p(\Gamma))$
is infinite dimensional. 
This shows the proposition.\qed

\bigskip

Following \cite{MS05}, we denote by
${\cal C}_{\rm geom}$ the class of countable groups
which admit a non-elementary weakly
acylindrical isometric action on some hyperbolic metric
space. Examples of such groups include.
\begin{enumerate}
\item[-] Word hyperbolic groups which are not virtually abelian.
\item[-] Any subgroup of the mapping class group of
an oriented surface of finite type and negative Euler characteristic
not preserving any essential multicurve, e.g. the
Torelli group.
\item[-] Any countable group which admits a non-elementary
isometric action on a (not necessarily locally finite) tree
which is proper on the edges.
\end{enumerate}
Our class also contains a large family
of \emph{relatively hyperbolic groups};
in fact, it seems that all geometrically finite
relatively hyperbolic groups in the sense
of Bowditch (see \cite{Y04} for a detailed 
discussion of those groups) are contained in ${\cal C}_{\rm geom}$.

For a locally compact $\sigma$-compact topological group $G$
define a \emph{lattice} in $G$ to be a \emph{discrete}
subgroup $\Gamma$ of $G$ such that $G/\Gamma$ admits
a \emph{finite} $G$-invariant measure. If $G=G_1\times G_2$
is any nontrivial direct product with locally compact
$\sigma$-compact and non-compact factors then we call
a lattice $\Gamma$ in $G$ \emph{irreducible} if
the projection of $\Gamma$ into each of the factors is dense.
The following lemma is part vi) of Proposition 7.13 
in \cite{MS05} and follows from the 
work of Burger and Monod \cite{BM02}.

\bigskip

{\bf Lemma 5.2:} {\it Let $\Gamma$ be an irreducible
lattice in a product $G=G_1\times G_2$ of locally compact
$\sigma$-compact non-compact groups; then $H_b^2(\Gamma,
\ell^2(\Gamma))=0$.}

\bigskip

We use Lemma 5.2 and the results of Monod and
Shalom \cite{MS05} to show.

\bigskip

{\bf Corollary 5.3:} {\it A group 
$\Gamma\in {\cal C}_{\rm geom}$ is not measure equivalent
to any finitely generated irreducible
lattice in either a simple Lie group of higher rank
or in a product of two locally compact $\sigma$-compact
and non-compact topological groups.}

{\it Proof:} By Theorem 4.5, for
every $\Gamma\in {\cal C}_{\rm geom}$ the group
$H_b^2(\Gamma,\ell^2(\Gamma))$ is nontrivial. Then 
Corollary 7.8 of \cite{MS05} shows that
$H_b^2(\Lambda,\ell^2(\Lambda))\not=\{0\}$ for 
every countable group $\Lambda$ which is measure
equivalent to $\Gamma$.

Now by Lemma 5.2, if $\Lambda$ is an irreducible lattice
in a product $G_1\times G_2$ of locally compact
$\sigma$-compact non-compact groups then
$H_b^2(\Lambda,\ell^2(\Lambda))=\{0\}$.
If $\Lambda$ is a lattice in a simple Lie group of
non-compact type and higher rank then the vanishing of the
second bounded cohomology group
$H_b^2(\Lambda,\ell^2(\Lambda))$ is due to Monod
and Shalom (Theorem 1.4 in \cite{MS}).
Thus in both cases, the group $\Lambda$ is not
measure equivalent to $\Gamma$.
(Note however that for lattices $\Lambda$ in simple Lie groups of higher
rank a much stronger result is due to Furman
\cite{Fu99a,Fu99b}: Every countable
group which is measure equivalent to $\Lambda$
is commensurable to $\Lambda$.)
\qed

\bigskip

Corollary C from the introduction now is immediate from
Corollary 5.3 and Proposition 5.1.

We conclude the paper with two additional applications
of the work of Monod and Shalom \cite{MS05}.

\bigskip

{\bf Corollary 5.4:} {\it A countable group containing
an infinite amenable normal subgroup is not measure
equivalent to a group in ${\cal C}_{\rm geom}$.}

\bigskip

Another consequence is Monod and Shalom's striking
rigidity result for actions of products
(Theorem 1.8 of \cite{MS05}).

\bigskip

{\bf Corollary 5.5:} {\it Let $\Gamma_1,\Gamma_2$ be torsion 
free groups in ${\cal C}_{\rm geom}$, $\Gamma=\Gamma_1\times 
\Gamma_2$ and let $(X,\mu)$ be an irreducible 
probability $\Gamma$-space. Let $\Lambda$ be
any torsion free countable group and let $(Y,\nu)$
be any mildly mixing probability $\Lambda$-space. 
If the $\Gamma$-action and the $\Lambda$-action are
orbit equivalent, then both groups as well as the
actions are commensurable.} 

\bigskip

There is also a version of Theorem A for
closed groups of isometries of proper hyperbolic
spaces and their 
continuous bounded cohomology \cite{H05b}.

\bigskip

{\bf Acknowledgement:} I am very grateful to Yehuda Shalom
for pointing out an error in an earlier version of
this paper and for additional valuable suggestions.
I also thank the anonymus referee for useful comments.

\bigskip

\noindent
MATHEMATISCHES INSTITUT DER UNIVERSIT\"AT BONN\\
BERINGSTRA\SS{}E 1 
\\53115 BONN, GERMANY

\smallskip

\noindent
{\it e-mail address}: ursula@math.uni-bonn.de

\end{document}